\documentclass[12pt]{article}
\begin{document}

\begin{titlepage}

\vspace*{-2cm}

\vspace{.5cm}

\begin{centering}

\huge{Chain extensions of $D$-algebras and their applications }

\vspace{.5cm}

\large  {Jining Gao }\\

\vspace{.5cm}

Department of Mathematics, North Carolina State University,
Raleigh, NC 27695-8205.

\vspace{.5cm}

\begin{abstract}

In order to unify the methods which have been applied to various topics
such as BRST
theory of constraints, Poisson brackets of local functionals, and certain
developments in deformation theory, we formulate a new concept which we
call the
{\it chain extension} of  a $D$-algebra.
We develop those aspects of this new idea which are central to
applications to algebra
and physics. Chain extensions may be regarded as  generalizations
of ordinary algebraic extensions of Lie algebras. Applications of our
theory provide
a new constructive approach to  BRST theories
which only contains three terms; in particular, this provides a new
point of view
concerning consistent deformations. Finally, we show how Lie algebra
deformations
are encoded into
the structure maps of an sh-Lie algebra with three terms. \end{abstract}

\end{centering}

\end{titlepage}

\pagebreak

\def\lh{\hbox to 15pt{\vbox{\vskip 6pt\hrule width 6.5pt height 1pt}
  \kern -4.0pt\vrule height 8pt width 1pt\hfil}}
\def\blob{\mbox{$\;\Box$}}
\def\qed{\hbox{${\vcenter{\vbox{\hrule height 0.4pt\hbox{\vrule width
0.4pt height 6pt \kern5pt\vrule width 0.4pt}\hrule height 0.4pt}}}$}}

\newtheorem{theorem}{Theorem}
\newtheorem{lemma}[theorem]{Lemma}
\newtheorem{definition}[theorem]{Definition}
\newtheorem{corollary}[theorem]{Corollary}
\newtheorem{proposition}[theorem]{Proposition}
\newcommand{\proof}{\bf Proof.\rm}

\section{Introduction}

Homological algebra has become an indispensable tool for the rigorous
formulation of
a wide variety  of developments in theoretical physics. Applications of
these
techniques
to physics has become so pervasive that  they have gradually become
identified as a new
category of mathematical physics which has been called ``cohomological
physics". One
of the fruitful branches of this theory is the ``cohomology"
formulation of the
BRST theory of constraints. Indeed the point of BRST theory is to
replace the
cohomology of the reduced space of a physical theory by the cohomology of a
homological resolution of the space being constrained. A separate
development
found in  \cite{BFLS}  applies homological techniques to classical
Lagrangian field as
a method of encoding Poisson brackets of local functionals.  Many of these
homological theories can be constructed explicitly when one is also
given the
existence of
a contracting homotopy.

Because of the diversity of approaches to these theories it is an attractive
challenge to find a more general algebraic structure to unify  existing
theories with the goal of obtaining even  more powerful applications of
these
theories. For example, in classic BRST theory, it is known how to
construct the BRST operator theoretically  via homological perturbation
theory(HPT),
but
explicit constructions of the BRST operator are difficult to obtain
because HPT
doesn't
tell us where the homological perturbation expansion terminates.
However, as one of the
applications of chain extensions
of $D$-algebras, this problem is settled in the case of  the Hamiltonian
formalism with
irreducible constraints, i.e., we have obtained  an explicit formula for
Hamiltonian
BRST  in
the most compact form yet known.
Another application in our paper is to provide a new investigation of
certain deformation problems such as consistent deformations and Lie algebra
deformations.
For example, in order to remove obstructions to Lie algebra deformation,
we embed
the Lie algebra
into an appropriate sh-Lie algebra  in such a way that the obstructions
will vanish
in the category of sh-Lie algebra deformations. Similarly, we can overcome
obstructions to consistent deformation
by embedding the original field into an appropriate extended superfield.
We expect these ideas  to have further applications in both physics and
algebra.

Our paper is organized as follows:
first of all, in section 2 we introduce the new and fundamental  concept
which we
call the chain extension of a $D$- algebra. In this section we also
develop the related
algebraic formalism needed for applications of these new constructs. In
section 3.1
we show how our new constructs provide a new approach to the BRST method of
encoding constraints on a symplectic manifold. In section  3.2 we show
how our new
methods apply to the homological approach to  local functionals in
Lagrangian field
theory referred to above. In section 3.3 we show how the consistent
deformation theory of Barnich and Henneaux relates to the chain extension
formalism.  In section 3.4 we show how certain  deformations of Lie
algebras are related to sh-Lie algebras. Finally, in a concluding
section we indicate
possible future development
of these ideas.

\section{Chain extension on $D$-algebras}

\begin{definition}

Let $V$ be a linear space over a field $k$ and $D$ a $k$-linear map such
that
$D^2=0$ on $V$. We call the ordered pair $(V,D)$  a
$D$-algebra.

\end{definition}

First we consider some important examples of $D$-algebras.
\bigskip

\noindent {\bf Example 1}:(Chain complex) Let $V$ be a graded space such
that
$V=\oplus_{p=0}^{\infty}V_p$ and $(D)_{p}$, $p\geq{0},$   a
sequence of maps (boundary operators)

$$\cdots V_{p+1} \stackrel{D_{p+1}}{\longrightarrow}V_{p} \cdots
\stackrel{D_2}{\longrightarrow}V_1 \stackrel{D_1}{\longrightarrow} V_0
\longrightarrow 0,$$
i.e.,  $D_p\circ D_{p+1}=0.$  For each $p,$ extend $D_p$ to the direct
sum $V$ by requiring
that it be linear and that $D_p$ restricted to $V_q$ be zero for $q\neq p.$
Define $D$ to be the direct sum $D=D_1+D_2+D_3+\cdots.$ The pair $(V,D)$
is then a
$D$-algebra.
\bigskip

\noindent {\bf Example 2}:(Lie algebra) Let $A$ be a linear space over
$k$, and let $l_2$
be a skew-linear map:
$l_2:A \otimes A\longrightarrow A$  which satisfies the  Jacobi identity:
\begin{eqnarray}
\sum_{\sigma \in unsh(2,1)}(-1)^{\sigma}  {l_2}( l_2(x_{\sigma(1)}\otimes
x_{\sigma(2)})\otimes x_{\sigma(3)})=0
\end{eqnarray}
where $x_1,x_2,x_3\in A.$ Then $l_2$ is a Lie bracket on $A;$ moreover
if  the left
side of above equation is  written in an abbreviated form  as  $$(l_2 l_2)(
x_1\otimes x_2\otimes x_3),$$
then the Jacobi identity is equivalent to the statement ${l_2}^2=0.$ Let
$V=\oplus_{n=0}^{\infty}{\otimes}A$
and extend the mapping $l_2:{A} \otimes {A}\longrightarrow {A}$ to a
mapping
$l_2:\otimes^{2+k}A\longrightarrow\otimes^{k}A$
via the equation

\begin{eqnarray}
 l_2(x_1\otimes \dots \otimes x_{2+k})= \nonumber \\ 
\sum_{unsh{(n,k)}}(-1)^{\sigma}e(\sigma) l_2 \left( x_{\sigma(1)}
\otimes x_{\sigma(2)} \right) \otimes
x_{\sigma(3)}\otimes\dots\otimes x_{\sigma(2+k)}.
\end{eqnarray}

Consequently we obtain an extended map $l_{2}:V \rightarrow V$ such that
${l_2}^2=0$
and
the pair  $(V,l_2)$ is a $D$-algebra.
\bigskip

\noindent {\bf Example 3}: Let $V$ be a linear space over $k$,
$TV=\oplus_{n=0}^{\infty}(\otimes^{n}V)$, and $D$  a coderivation
on the cofree coalgebra $TV$ which satisfies $D^2=0$. Then $(TV,D)$ is a
$D$-algebra. All $A_{\infty}$ algebras are of this type. When $V$ is a
graded space
one may obtain the class of all sh-Lie algebras by selecting a subspace
of the space
$TV$ and by modifying the coalgebra structure on the subspace so that
the sh-Lie
structure is characterized by a coderivation $D$ with square zero. Thus
they are
also included as a special type of $D$-algebra.
\bigskip

Let $({\cal F},D_{\cal F})$   be a $D$-algebra over a field $k$ and
${\cal B}$ be an
arbitrary linear space over $k$.
Consider the ordered triple $({\cal B},{\cal F},D_{\cal F})$. The
primary new
concept in
this paper is the notion of a chain extension of such a triple.

\begin{definition}

Let $({\cal B}, {\cal F},D_{\cal F})$   be a triple defined as above and let
$(X_*,\delta)$ be a homological resolution of ${\cal F}$, i.e.,  we have
a complex
$$\cdots X_{p+1} \stackrel{\delta_{p+1}}{\longrightarrow}X_{p} \cdots
\stackrel{\delta_1}{\longrightarrow}X_1
\stackrel{\delta_0}{\longrightarrow} X_0 $$
where $X_0={\cal B}\oplus{\cal F},  \delta _{0}(X_1)={\cal B} ,$
and  $H_*(X)=H_0(X){\simeq}{\cal F}$.
If there exists a sequence of maps $l_n :X_*{\rightarrow}X_*$ with degree
$n-2$ such that $l_1=\delta,  l_2{\mid}_ {\cal F}=D_{\cal
F} ,$ and such that  the summation $ l=\sum_{n=0}^{\infty}l_n$ is a
nilpotent operator on
$X_*,$ then we call $(X_*,l)$ a chain
extension of ${\cal B}$  by $({\cal F},D_{\cal F})$ .
\end{definition}

Notice that our resolution
$\{X_n\}$ is defined only for $n\geq 0,$ thus we do not consider
augmented complexes
in this paper. Although the definition above initially appears quite
complicated, we
will find that
the notion of a chain extension is a generalization of ordinary
algebraic extension theory. Indeed ordinary algebraic extensions can be
described as
follows.
Consider a short exact sequence:
\begin{eqnarray}
& & 0\longrightarrow C\stackrel{i}{\longrightarrow}B
\stackrel{j}{\longrightarrow}A\longrightarrow 0
\end{eqnarray}
where $C,B,A$ are linear spaces over $k$ and where one has a
$D$-algebra structure
on $B.$  Furthermore assume that the restriction
of $D$ on $C$ defines a $D$-algebra structure on $C$ and induces a
mapping $\bar D:
A\rightarrow A$ such that  $A$ is a $\bar D$-algebra. Also notice that
the chain
complex
$0\longrightarrow C\stackrel{i}{\longrightarrow}B\longrightarrow 0$
provides a
resolution of $A$
for which $\delta_k=0$ for $k>0,$ $\delta_0=i.$ If the quotient mapping
$j:B\rightarrow A$ plays the role of $\eta$ in the diagram below,
then $S=l_1+l_2=\delta+D$ is a chain extension of $C$ by $A.$ Such
$D$-algebra
extensions of $C$ by $A$ occur in some applications in physics.

Sometimes a homological resolution for $\cal F$ comes with the
additional structure
of a contracting
homotopy which we will show often guarantees the existence of a chain
extension.  Furthermore the chain extension in this case may be defined
explicitly
in terms of this homotopy. We adapt the constructions in
\cite{BFLS} to describe such a contracting homotopy.

Note that  $\cal F$ may be regarded as a  differential graded vector
space ${\cal
F}_*$ with ${\cal F}_0={\cal F}$,  and ${\cal F}_k=0$ for
$k>0.$ Assume that there exists a chain map $\eta:X_*\longrightarrow
{\cal F}_*$
with homotopy inverse
$\lambda:{\cal F}_*\longrightarrow X_*$;  i.e., we have that $\eta \circ
\lambda =
1_{{\cal F}_*}$ and that $\lambda \circ \eta
\sim 1_{X_*}.$  Thus there is  a chain homotopy $s:X_*\longrightarrow
X_*$ with
$\lambda \circ \eta -1_{X_*}=l_1\circ s +s \circ l_1$. In particular,
notice that in
degree zero one has that $1_{X_*}=\lambda \circ \eta - l_1\circ s.$

We may summarize all of the above with the commutative diagram

$$
\begin{array}{cccccc}
&&s& &s\\
\cdots\longrightarrow & X_2 &\stackrel
{\textstyle\longleftarrow}{\longrightarrow} &X_1&\stackrel
{\textstyle\longleftarrow}{\longrightarrow} &X_0\\ & & l_1 & & l_1&\\
& \lambda\Bigl\uparrow\Bigr\downarrow\eta & &
\lambda\Bigl\uparrow\Bigr \downarrow\eta&&
\lambda\Bigl\uparrow\Bigr\downarrow\eta\\ &&&&&\\
\cdots\longrightarrow &0&\longrightarrow & 0 &\longrightarrow &
H_0={\cal F}.
\end{array}
$$

\begin{theorem} Assume that $({\cal F},D_{\cal F})$ is a $D$-algebra and
that the pair
$({\cal F},D_{\cal F})$ has a homological resolution with contracting
homotopy $\eta$ as
indicated above. If  $ {\cal B}=l_1(X_1) $ and there exists a mapping
$l_2:X_0
\rightarrow X_0$ which satisfies the properties:
\begin{eqnarray}
&&   (i) \quad  l_2|_{\cal F}=D_{\cal F}\\
& &  (ii)\quad  l_2(\cal B)\subseteq  \cal B \label{g}\\
& & (iii)\quad  l_{2}^{2} (X_0)\subseteq \cal B \label{h}.
\end{eqnarray}
then there exists a chain extension  $(X_*,l)$ of the triple $({\cal
B},{\cal
F},D_{\cal F})$ which  can be  explicitly constructed by induction as
follows. For
any $x \in X_*:$\newline
(1) define $l _1(x)=\delta (x) ,$
\newline
(2) for deg $ x>0,$ define $l_2(x)=(s\circ (l_2l_1))(x),$  \newline
(3) if deg $x=0$, define $l_3(x)=(s\circ (l_2 l_2))(x)$ and for deg $x>0$
define $l_3(x)=(s\circ (l_2 \circ l_2+l_3 \circ  l_1))(x),$
\newline
(4) for all $x$ and for $n>3,$ define $l_n(x)=0.$\newline
Using these definitions it follows  that  $ l_2(x)=0$ for deg $x>1$ and
that
$l_3(x)=0$ for deg $x>0.$
\end{theorem}

${\bf Remark :}$ Conditions $(i)$ and $(ii)$ of the theorem above were
first introduced in \cite{BFLS} where
an sh-Lie prototype of the theorem above was proved. Markl  was first to
notice (via
a private communication) that  these conditions insure that at most
three of the
sh-Lie structure maps are nontrivial. Barnich proved this result in
\cite{B}.
Finally, Al-Ashhab \cite{Samer} weakened the conditions of Markl's
original statement and provided the
details of  his observations to the author who generalized the result to
the case of
$D$-algebras, for which we have more extensive applications.

\begin{proof}
Notice first that because ${\cal B}=l_1(X_1),$ it follows from the
commutative diagram
above that $\eta({\cal B})=0.$ This fact will be used throughout the
proof.  Assume
that there exists an operator $l_2:X_0\rightarrow X_0$ which satisfies the
hypothesis
of the theorem.

In order to prove  the operator  $l=l_1+l_2+l_3$ is nilpotent, observe
that if we
expand
$l^2$ the equality $l^2=0$ is equivalent to the sequence of
equations.
\begin{eqnarray}
l_1l_2+l_1l_2=0\label{b}\\ {l_2}^2+l_1l_3+l_3l_1=0 \label {c} \\
l_2l_3+l_3l_2=0
\label {k}
\\ {l_3}^2=0  \label{d}\\ \nonumber
\end{eqnarray}
We give an inductive proof of these equations. For any $x\in {X_0}$, since
$l_1l_2(x)=l_2l_1(x)=0$, we have
$(l_1l_2+l_1l_2)(x)=0$, and  (\ref{b}) holds in degree zero.

We now prove (\ref{c}) in degree zero, recall that  $\lambda \circ \eta
-1_{X_*}=l_1\circ s +s \circ l_1$, and apply both sides of this equation
to $({l_2}^2+l_3l_1)(x)$ ,where $x\in {X_0}.$
Now $l_1(X_0)=0$ and $l_2^2(X_0)\subseteq {\cal B}$ implies that,
$(\lambda \circ \eta )({l_2}^2(x))=0,$
consequently
$$(\lambda \circ \eta -1_{X_*})({l_2}^2+l_3l_1)(x)(\lambda \circ \eta )(
{l_2^2}(x))-l_2^2(x)=-l_2^2(x).$$
For similar reasons,  $(s\circ l_1)({l_2}^2+l_3l_1)(x)=0$ and the
equation above
implies
\begin{eqnarray}-l_2^2(x)=(l_1\circ s +s \circ
l_1)({l_2}^2+l_3l_1)(x)\\=(l_1\circ
s)({l_2}^2+l_3l_1)(x)
=l_1(s(l_2^2(x)))=l_1(l_3(x)),
\end{eqnarray}
where in the last equality we used the inductive definition of $l_3.$
This shows that $(-l_2^2)(x)=(l_1l_3)(x)$ and since $(l_3l_1)(x)=0,$
we have \begin{eqnarray}({l_2}^2+l_1l_3+l_3l_1)(x)=0 \end{eqnarray}
Thus (\ref{c}) holds for $x \in {X_0} $.

We now give an inductive proof of (\ref{b}) in arbitary degree. Assume
that (\ref{b})
holds for any $x'\in X_*$ with deg $x'=n$;  if deg $x=n+1$, we
first prove that $(l_2\circ l_1)(x)$
is an $l_1$ boundary.  Since deg $l_1(x)=n,$ it follows from  (\ref{b})
that
\begin{eqnarray}
l_1 (l_{2 } l_1)(x)=({l_1}  l_2)(l_1(x))=(-l_2{l_1})(l_1(x))=0
\end{eqnarray}
Thus $(l_2l_1)(x)$ is $l_1$-closed. Notice that by the inductive
definition $l_2$
preserves
degree and so deg($(l_2l_1)(x))=n.$ If $n>0,$  then the cycle
$(l_2l_1)(x)$ is an
$l_1$-boundary since $H_n(X_*)=0.$ If $n=0$ then deg $x=1,  l_1(x)\in
{\cal B},$ and
$(l_2l_1)(x)\in l_2({\cal B})\subseteq {\cal B}=l_1(X_1).$ Thus
$(l_2l_1)(x)$ is
also a boundary in this case.
Now apply $\lambda \circ \eta -1_{X_*}=l_1\circ s +s \circ l_1$ to
$l_2{l_1}(x);$ we
have
\begin{eqnarray}
& & (\lambda \eta)(l_2 l_1)(x)-l_2 l_1(x)=l_1 s(l_2
l_1(x))+sl_1(l_2{l_1}(x))
\end{eqnarray}
Since $(l_2{l_1})(x)$ is a boundary and $\eta (l_2{l_1})(x)=0$,
\begin{eqnarray}
l_2{l_1}(x)+l_1{s}l_2{l_1}(x)=0.
\end{eqnarray}
Since deg $x>n+1,$ it follows from the inductive definition of $l_2$ that
$s(l_2l_1)(x)=l_2(x)$ and therefore that $l_2l_1(x)+l_1l_2(x)=0.$
Thus (\ref{b}) holds.

In order to  prove (\ref{c}) for arbitrary degree, we also need to first
verify that
for any $x\in {X_*}$,
$l_2^{2}(x)+l_3{l_1}(x)$ is a boundary.
Since $l_2^{2}(X_0)\in {\cal B}$, we immediately get that
$l_2^2(x)=l_2^{2}(x)+l_3{l_1}(x)$ is a boundary
for all $x\in {X_0}.$  Assume that it is a boundary  for every $x'$ with deg
$x'<n+1$. Choose any $x$ such that deg $x=n+1$ and observe that
\begin{eqnarray}
l_1(l_2^{2}+l_3{l_1})(x) & = & (l_1{l_2^{2}})(x) + (l_1{l_3}{l_1})(x) \nonumber\\
 & = & (l_1{l_2^{2}})(x) + (l_1{l_3})(l_1(x))
\end{eqnarray}

Since deg $l_{1}(x)=n$, it follows from  the inductive hypothesis
applied to (\ref{c})  that
$l_{1}l_{3}l_{1}(x)=-(l_{2}^{2}+l_{3}l_{1})l_{1}(x)=-l_{2}^{2}l_{1}(x);$
thus $l_{1}(l_{2}^{2}+l_3l_1)(x)=(l_{1}l_{2}^{2})(x)-l_2^{2}l_1(x)=0$
(since it
follows from (\ref {b})
that $l_1$ and $l_2$ anti-commute).
It follows that $l_2^{2}(x)+l_3l_1(x)$ is a cycle. It also follows from the
inductive definition of $l_3$ that
it increases the degree by 1, thus deg ($l_2^{2}(x)+l_3l_1(x))=n+1.$ Since
$H_{n+1}(X_*)=0,$ $l_2^{2}(x)+l_3l_1(x)$ is a boundary. Applying
$\lambda \circ
\eta -1_{X_*}=l_1\circ s +s \circ l_1 $ to $l_2^{2}(x)+l_3l_1(x),$ we have
\begin{eqnarray}
(\lambda
{\eta})(l_{2}^{2}+l_3{l_1})(x)-(l_{2}^{2}+l_3{l_1})(x)=l_1{s}(l_2^{2}(x)+l_3{l_1}(x))+sl_1(l_2^{2}(x)+l_3{l_1}(x))\nonumber
\end{eqnarray}
and $l_{2}^{2}(x)+l_3l_1(x)+l_1{s}(l_2^{2}+l_3l_1)(x)=0.$ By the inductive
definition of $l_3,$
$s((l_2^2+l_3l_1)(x))=l_3(x),$ and consequently
$({l_2}^2+l_3l_1+l_1l_3)(x)=0.$  Thus (\ref{c}) follows.

We now prove the two remarks in the last statement of the theorem.
Assume that deg
$x>1,$ then  deg $x\geq 2$ and deg $l_1(x)>0.$, consequently
$l_2(x)=sl_2l_1(x)=s(sl_2l_1)(l_1(x))=0.$
The first of the two statements follows.  Now assume that deg $x\geq 1$,
then using
the inductive definitions of both $l_2$ and $l_3,$ we have
\begin{eqnarray}
l_3(x)=s(l_2^{2}+l_3{l_1})(x)\nonumber\\=
s((sl_2{l_1})l_2+s(l_{2}^{2}+l_3{l_1})l_1)(x)\nonumber\\
=s^2((l_2{l_1})l_2+(l_{2}^{2}l_1+l_3{l_1}^2))\nonumber\\
=s^2[-l_1l_2^{2}+l_2^{2}l_1(x)]=0  \nonumber
\end{eqnarray} and the second statement follows.
It is now easy to prove (\ref{k}),(\ref{d}) as follows:
since $l_3(x)=0 $ for any deg $x>0,$  it's enough to prove
(\ref{k}),(\ref{d}) in degree
zero.
Assume that $x\in {X_0}$, we have deg $l_3(x)$=1, thus $l_3l_3(x)=0 ,$
and (\ref{d})
follows.

We now prove (\ref{k}): let $x\in {X_0},$ since deg $l_2(x)=0,$ deg
$l_3(x)=1,$ and by
the inductive definitions of $l_2,l_3 ,$ we have

$l_2l_3(x)+l_3l_2(x)=sl_2l_1l_3(x)+sl_2l_2l_2(x)\\
=sl_2(l_1l_3+l_2l_2)(x)=sl_2((l_1l_3+l_2l_2+l_3l_1)(x)=0.$

This concludes the proof that $l^2=0$ and the fact that $(X_*,l)$ is a
chain extension.
\end{proof}

From the theorem above, we find that a chain extension for a given
triple $(\cal B,
\cal F, D)$ is not unique,but the following theorem
will tell us all chain extensions have the same homology  group.
\begin{theorem}
Let $(X_*,l)$ be a chain extension of triple $(\cal B, \cal F, D)$, then
$H(X_*,l)\simeq H(\cal F,D)$.
\end{theorem}

The proof of this theorem is very similar to the homological
perturbation technique
of the theorem on page 179-180 in \cite{HT}.

\section{Applications}
\subsection{BRST theory}
As one of the applications of  the theorem above, we provide a new
construction for the
Hamiltonian BRST operator  in the irreducible case, which has the most
compact form
possible(three maps are enough!).  In \cite{HT},  a $BRST$ operator is
obtained by "homological perturbation theory" but not explicitly. Using
our Theorem 3,
we will  reconstruct a $BRST$
operator in an explicit way under some local conditions. In particular,
we show that  the BRST theory for the irreducible case  is just a chain
extension of the
longitudinal differential.

Let $P$ be a symplectic manifold, $G_a=0$,$(a=1,\cdots ,{n})$ a set of
first class
constraints defined on $P,$ and $\Sigma$  the constraint surface defined
by the zeros of the functions  $\{G_a\}$. We assume that $\Sigma$ is a
regular submanifold of
$P.$ Let  $N$ be the ideal of the algebra $C^{\infty}(P)$ generated
by the constraint functions $G_a, a=1,\cdots,{n}$ and consider the short
exact sequence:
\begin{eqnarray}
0\longrightarrow{N}\stackrel{i}{\longrightarrow}
C^{\infty}(P)\stackrel{j}{\longrightarrow}C^{\infty}(\Sigma)\longrightarrow
0. \label{i}
\end{eqnarray}

In general, the sequence (\ref{i}) is not split as an algebra over the
complex field $C,$  but since $\Sigma$ is a regular
submanifold there exists a local chart in $P$
on which the sequence is split locally. That is to say,  if $P$ is
replaced by a
suitable open subset of $P,$ there exists a splitting algebra homomorphism
$\tau:C^{\infty}(\Sigma)\rightarrow C^{\infty}(P)$ such that one obtains
the direct
sum decomposition $C^{\infty}(P)=N\oplus C^{\infty}(\Sigma)$ as algebras.

If $[\cdot,\cdot]$ is the Poisson bracket defined by
the symplectic structure  on $P,$ then we may define a mapping on
$C^{\infty}(\Sigma)=C^{\infty}(P)/N$ by $g+N \rightarrow [g,G_a]+N$ for
each $a.$
This mapping is a well-defined derivation of $C^{\infty}(\Sigma).$ We
denote its value
at $f\in C^{\infty}(\Sigma)$ by $\partial_af$ but sometimes by an
obvious abuse
of notation we denote it by $[f,G_a].$

Let ${\cal F}=C^{\infty}(\Sigma)\otimes C(\eta^{a})$ where  $C(\eta^a)$
is the exterior algebra
over the complex field $C$ with generators   $\{\eta^a\}$ of degree
one.  Physicist usually refer
to  the generators  $\eta^a$ as ghost variables.

Define the longitudinal derivative $d$ on ${\cal F}$  by defining it on
generators of ${\cal F}$ and then
extending it to all of ${\cal F}$ as a right derivation. On generators it  
is
defined as follows(\cite{HT}):
\begin{eqnarray}
df=[f,G_a]\eta^ {a} \label{diff1}\\
d\eta^a={\frac{1}{2}}C_{cb}^{a}\eta^{b}\eta^{c}
\label{diff2} \end{eqnarray}
where $f \in C^{\infty}(\Sigma)$ and  $C_{cb}^{a}$ are structure
functions on $P.$

Notice the defining equations of $d$ on generators given above is
induced from the action of $d$
on $C^{\infty}(P)\otimes C(\eta^a)$ given by the same equations (due to
our abuse of the use of the
bracket referred to above). We claim that
$d^2(C^{\infty}(P)\otimes C(\eta^a))\subseteq N\otimes C(\eta^a).$ The
proof of this fact follows from induction on the ghost number. For the
convenience of the reader
we show that this is true for ghost number one. The inductive step is
straightforward and is left
to the reader. Let $\alpha=f_b\eta^b$ where
$f_b\in C^{\infty}(P),$ and observe that
$d^2\alpha=(d^2f_b)\eta^b+f_bd^2\eta^b.$
Now for $f\in C^{\infty}(P),$
$$ d^2f=- [[f,G_a],G_b]\eta^b\eta^a -\frac
{1}{2}C^a_{bc}[f,G_a]\eta^b\eta^c.$$
It follows from the Jacob identity and the anti-commutativity  of the
$\eta$'s that
$$2[[f,G_a],G_b]\eta^b\eta^a=[f,C^c_{ab}G_c]\eta^b\eta^a.$$
Thus
$d^2f=-\frac{1}{2}[f,C^c_{ab}]G_c\eta^b\eta^a-\frac{1}{2}C^c_{ab}[f,G_c]\eta^b\eta^a
-\frac {1}{2}C^a_{bc}[f,G_a]\eta^b\eta^c$
from which it follows that
$d^2f=-\frac{1}{2}[f,C^c_{ab}]G_c\eta^b\eta^a,$ which is obviously
in $ N\otimes C(\eta^a).$

Thus $d^2f$ belongs to $N\otimes C(\eta^a)$ for each $f\in
C^{\infty}(P)$ as does also $(d^2f_b)\eta^b$ above. Recall that when the
$C^a_{bc}$ are constant one  has
that $d^2\eta^a=0.$ On the other hand when the structure constants are
actually structure
functions, then further calculations similar to those above using
(\ref{diff1}),(\ref{diff2}) and
anti-commutivity of the $\eta$'s  show that
$$d^2\eta^c=\frac{1}{2}\{ C^c_{ab}
C^b_{pq}+[C^c_{aq},G_p]\}\eta^a\eta^p\eta^q.$$
Moreover a similar calculation shows that
$$[[G_a,G_p],G_q]\eta^a\eta^p\eta^q=-\{C^c_{ab}
C^b_{pq}+[C^c_{aq},G_p]\}G_c\eta^a\eta^p\eta^q.$$
It follows that $-2d^2\eta^cG_c=[[G_a,G_p],G_q]\eta^a\eta^p\eta^q.$ The
latter is zero by the
Jacobi identity and the anti-commutativity of the $\eta$'s.
Since $d^2\eta^cG_c=0,$ it follows from the fact that the constraints
are irreducible that $d^2\eta^c$ is in the ideal $N\otimes C(\eta^a)$
for each $c.$

It follows that $d^2\alpha \in N\otimes C(\eta^a)$ for each $\alpha\in
C^{\infty}(P)\otimes C(\eta^a)$
of ghost degree 1. The general case now follows easily by induction on
the ghost degree.

From this it follows that $({\cal F},d)$ is a $D$-algebra. We intend to
show that the $BRST$ complex can be regarded
as a chain extension of this algebra. In order to construct the
homological resolution of $\cal F,$ we need an antighost variable $P_a$
for each
constraint $G_a$ and a Kozul-Tate differential $\delta$
  defined as follows:
\begin{eqnarray}
\delta {f}=0\\{\delta}{P_a}=-G_a\\{\delta} {\eta^ {a}=0}
\end{eqnarray}
where $f\in C^{\infty}(P).$

Let $X_*=C(P_a) \otimes C^{\infty}(P) \otimes
C(\eta^{a})$ be the graded space where the grading
is given by the antighost degree and let ${\cal F}_*$ be the graded
space with
${\cal F}_0={\cal F}$ and with ${\cal F}_k=0, k>0.$
If we extend the longitudinal derivative d to the anti-ghosts
by requiring that $dP_a=0$ for all $a,$ then we can extend both $d$ and
$\delta$ to all of
$X_*$ by requiring that both of them be right derivations. We will show that
$(X_*,\delta)$ is a resolution of $\cal F.$

Define ${\cal B}$ by  ${\cal B}=\delta(X_1)=N\otimes C(\eta^{a})$ and let
$\tilde \eta: X_*\rightarrow {\cal F}_*$ be the chain map which is
nontrivial only in degree
zero and in that case define it to be the projection
$$X_0\rightarrow X_0/{\cal B}=(C^{\infty}(P) \otimes C(\eta^{b}))/(N
\otimes C(\eta^{b}))
=C^{\infty}(\Sigma) \otimes C(\eta^{b})={\cal F}.$$

We now construct a contracting homotopy, but to do this we need to
modify the contracting homotopy formula in \cite{HT} to
conform to the conventions in our paper. To make contact with \cite{HT}
we adopt their
notation throughout the remainder of this section.  In particular
we employ  the formal algebraic language used by physicist in order to
show that their formulation
of Hamiltonian BRST theory may be reformulated in terms of chain
extensions.

Since
the constraints $G_a=0$ are assumed to be independent (recall that we
are working in
the irreducible case), we can choose a local chart
$(x_i,G_a)$ in a neighborhood $O$ of the  manifold P. Our results are
valid only on
such an
open set $O$ and since all the structures of interest restrict to $O$ we
presume
that $O=P$ for simplicity. Let
\begin{eqnarray}
\delta=(\frac{\partial^{R}}{\partial{P_a}})G_a ,& &
\sigma=(\frac{\partial^{R}}{\partial{G_a}})P_a
\end{eqnarray}

\begin{eqnarray}
\bar
N=(\frac{\partial^{R}}{\partial{P_a}})P_a+(\frac{\partial^{R}}{\partial{G_a}})G_a
\end{eqnarray}
where $\frac{\partial^{R}}{\partial{P_a}}$ and
$\frac{\partial^{R}}{\partial{G_a}}$
are right derivations on the algebra $ C(P_a)\otimes{C^{\infty}(P)}\otimes
{C(\eta^b)}.$

We will use the mapping $\sigma$ to build the contracting homotopy $s$
below. To
obtain the necessary properties of $s$ parities are important. In the
calculations below
we assume that the parities of $\delta, \sigma, \bar N$ are all zero.
Moreover we require
the parity of smooth functions on $P$ be zero. The
parities of ghost fields, $\varepsilon(\eta^a),$ and anti-ghost fields,
$\varepsilon(P_b),$
are related to the parities of functions on $P$ by $\varepsilon(\eta^a)=
\varepsilon(G_a)+1= \varepsilon(P_a).$
The relation $\delta{\sigma}+\sigma{\delta}=\bar{N}$ holds when
evaluated at an
arbitary element of $ C(P_a)\otimes{C^{\infty}(P)}\otimes {C(\eta^b)} $
and the
latter equation implies that:
\begin{eqnarray}
t\frac{d}{dt}F(tP_a,tG_a,x_i,\eta^b)=(\delta{\sigma}+\sigma{\delta})F(tP_a,tG_a,x_i,\eta^b)
\label{j}
\end{eqnarray}
where $t$ is an arbitary real number. If, in addition, F is of antighost
number
$k>0$, it follows from  (\ref{j}) that
\begin{eqnarray}
& & F(P_a,G_a,x_i,\eta^b)=(\delta \sigma+\sigma
\delta)(\int_{0}^{1}\frac{F(tP_a,tG_a,x_i,\eta^b)}{t} dt).
\end{eqnarray}
If we set $\psi(F)=-\int_{0}^{1}\frac{F(tP_a,tG_a,x_i,\eta^a)}{t} dt$  and
$s=\sigma\circ {\psi},$ then $\delta\circ \psi =\psi \circ \delta$ and
$l_1\circ s+s\circ l_1=-1_{X_*}$ in degree $k>0.$ Now consider the degree
zero case. Recall that in the definition of chain homotopy one is required 
to
find a homotopy inverse $\lambda:{\cal F}_*\rightarrow X_*$ of the chain 
mapping
$\tilde \eta: X_*\rightarrow {\cal F}_*.$  Both $\tilde \eta$ and $\lambda$ 
should
be   nontrivial only in degree
zero in which case $\tilde \eta$ is the projection
$X_0\rightarrow X_0/{\cal B}={\cal F}$ defined above. Thus we need to define
both $\lambda$ and $s$ in degree zero in such that
$\lambda \circ \tilde \eta -1_{X_*}=l_1\circ s +s \circ l_1.$

A degree zero element of $X_*$ is simply  a function $f$ such that
$f(G_a,x_i,\eta^b) \in C^{\infty}(P) \otimes C(\eta^{b}) .$ Define
$sf= \int_{0}^{1} \frac{1}{t} (\frac{\partial^R}{\partial G_a}) 
f(tG_a,x_i,\eta^b)P_a
dt.$
In order to construct the map $\lambda $ for the case of $antigh(f)=0,$ we 
first
consider
the map $\tilde {\lambda} : X_0 \longrightarrow X_0 $ defined as follows.
For any $f \in X_0,$ let
\begin{eqnarray}
& & \tilde {\lambda}(f)=f+\delta sf
\end{eqnarray}
To find the desired map the following Lemma is useful.

\begin{proposition}
$\tilde {\lambda} $ vanishes on the subspace ${\cal B}$ and thus induces a 
map
$\lambda: X_0 /{\cal B} \longrightarrow X_0 $
\end{proposition}
\begin{proof} Since $ f\in {\cal B},$ the integral
$sf=\int_0^1\frac{1}{t}(\frac{\partial^R}{\partial G_a}) 
f(tG_a,x_i,\eta^b)P_a dt$
exists
and we have:
\begin{eqnarray}
& & \tilde {\lambda}(f)=f+\delta sf \nonumber \\
& &=f(G_a,x_i,\eta^b) - \int_{0}^{1} \frac{1}{t}(\frac{\partial^R}{\partial 
G_a})
f(tG_a,x_i,\eta^b) G_a dt\nonumber \\
& & =f(G_a,x_i,\eta^b)-\int_{0}^{1}\frac{d}{dt}f(tG_a,x_i,\eta^b) dt 
\nonumber \\
& & =f(G_a,x_i,\eta^b)-(f(G_a,x_i,\eta^b)-f(0,x_i,\eta^b))\nonumber \\
& & =f(0,x_i,\eta^b)=0 \nonumber \end{eqnarray}
where the last step follows from the fact that  $f \in {\cal B} .$
\end{proof}

It is now easy to verify that the maps $s,\lambda, \eta $ satisfy the 
contracting
homotopy
formula on the space $X_0$.

To summarize, the contracting homotopy is obtained as follows: given any
$F(P_a,G_a,x_i,\eta^b)\in C(P_a)\otimes C^{\infty}(P) \otimes
C(\eta^{b}),
s(F)=\int_{0}^{1} \frac{1}{t}(\frac{\partial^R}{\partial G_a}) 
f(tG_a,x_i,\eta^b)P_a
dt$
when $ antigh(f)=0$  and
$s(F)=\sigma{\psi}(f)$ for $antigh(f)>0.$
It follows that we have a homological resolution $(X_*,\delta)$ with
contracting homotopy  $s.$

We now show that there exists maps $l_1,l_2,l_3$ which satisfy the
definition
of chain extension. To do this, we show that if  $l_1=\delta$ and $ l_2=d$
then the hypothesis of Theorem 3 is true. First observe that
$X_1=V(P_a)\otimes C^{\infty}(P)  \otimes C(\eta^{b})$ where $V(P_a)$ is
the linear
space spanned by the $P_a$'s over $C.$ Recall that
$\delta(X_1)=N\otimes C(\eta^{b})={\cal B}$
and observe
that $l_2({\cal B})=d(N \otimes C(\eta^{b}))\subseteq N\otimes
C(\eta^{b})={\cal B}.$
Moreover, note that $X_0=C^{\infty}(P)\otimes C(\eta^a)$ and recall
that  we have already
shown that $d^2(C^{\infty}(P)\otimes C(\eta^a))\subseteq N\otimes
C(\eta^a)$ in our proof
that $({\cal F},d)$ is a $D$-algebra.
Thus $l_2^2(X_0)\subseteq N \otimes C(\eta^{b})= {\cal B}.$ It follows
from Theorem 3
that there are maps $l_1,l_2,l_3$ which satisfy the definition of a
chain extension and that
they are defined inductively by the Theorem.

We can now use Theorem 3  to determine the   $BRST$ operator.
Let's start with the definition of
$l_2$ on $X_1;$ obviously knowing the values of $l_2$ on $X_1$ is enough
to determine $l_2$ on the whole space
$X_*$ because elements of degree greater than zero are  generated by
$X_1.$ By
applying formula (2) of Theorem 3, we have
\begin{eqnarray}
l_{2}P_a=dP_a=sl_2 l_1(P_a)=sl_2\delta(P_a)\nonumber\\
=sl_2(-G_a)\nonumber\\
=-s((\partial_b{G_a})\eta^b)\nonumber\\ =-s([G_a,G_b]\eta^b)\nonumber\\
=-s(C_{ab}^{c}G_{c}\eta^{b})\nonumber\\  =-s(C_{ab}^{c})G_c
\eta^{b}-C_{ab}^{c}s(G_c)\eta^{b}\nonumber\\
  =-s(C_{ab}^{c})G_c
\eta^{b} -C_{ab}^{d}P_{d}\eta^b
\end{eqnarray}
The next operator we need to compute is $l_3,$ and by Theorem 3, we
need only
compute it in the case that
$ {f=f(x_i,G_a,\eta^b)}$ is in $ C^{\infty}(P)\otimes C(\eta^b).$ Recall
that
$Y_a=\partial_a$ are the Hamiltonian vector fields of the constraints
$G_a$ and that
$[\partial_a,\partial_b]=C^c_{ab}\partial_c+Y_{C^c_{ab}}G_c.$ By  (3) of
Theorem 3,
we have:
\begin{eqnarray}
l_3(f)=s(l_2^2+l_3l_1)(f)\nonumber\\=s(l_2^2)(f)\nonumber\\=
sl_2((\partial_a{f})\eta^a)\nonumber\\
=s(l_2(\partial_a{f})\eta^{a}+(\partial_a{f})l_2(\eta^{a}))\nonumber\\
=s(\partial_b(\partial_{a}f)\eta^{b}\eta^{a}+(\partial_a{f})(-\frac{1}{2}C_{bc}^{a}\eta^{b}\eta^{c}))\nonumber\\
=s(\frac{1}{2}([ \partial_a, \partial_b]f) \eta^{b}
\eta^{a}-\frac{1}{2}(\partial_c{f}) \eta^{b} \eta^{a)} \nonumber\\
=s(\frac{1}{2}([\partial_a,\partial_b]f) -
C_{ba}^{c}\partial_c)f)\eta^b\eta^{a})\nonumber\\
=s(\frac{1}{2}(G_c(Y_{C_{ba}^{c}}f) \eta^b \eta^a)\nonumber\\
=\frac{1}{2}s(G_c)(Y_{C_{ba}^{c}}f)\eta^{b}
\eta^{a}+\frac{1}{2}G_cs(Y_{C_{ba}^{c}}f)\eta^{b} \eta^{a}\nonumber\\
=\frac{1}{2}(Y_{C_{ba}^{c}}f)P_c \eta^{b}
\eta^{a}+\frac{1}{2}G_cs(Y_{C_{ba}^{c}}f)\eta^{b} \eta^{a} \nonumber\\
=\frac{1}{2}[C_{ab}^{c},f]P_c\eta^{b}\eta^{a}+\frac{1}{2}G_cs(Y_{C_{ba}^{c}}f)\eta^{b}
\eta^{a} \nonumber\\
\end{eqnarray}
and $l_3(P_a)=l_3(\eta^{a})=0$

${\bf Remark :}$ From the observations above, we see that our new
development of the
BRST operator
has only three nonzero terms. Moreover, it is clear that this is the
least number of
terms required to describe the nontrivial case we have focussed on here,
i.e., in this case
we have acheived an optimal BRST construction. Our calculation of the
$BRST$ operator above
depends on the existence
of an explicit contracting homotopy.  However, the
contracting homotopy formula required by our construction only holds
locally,  in
other words, it
depends on the choice of local coordinates.
In the $BRST$ theory of Lagrangian field theory, such a  choice of a
local chart is
referred to as a "regularity condition". General speaking, conditions of
this type
impose restrictions on the constraint surfaces (functions). Further
details regarding regularity conditions can be found on pages 7 and 199
in \cite{HT}.
Such conditions are
given  explicitly for field theories  such as the  "Klein-Gordon"
and Dirac fields defined on appropriate spacetimes. We can apply the
method above to
examples of such field theories.

\subsection{Local problem via chain extensions}

The formulation of Lagrangian field theory in terms of jets of sections
of vector
bundles
has proven to be quite satisfactory in the modeling of large classes of
physical
theories when
the fields are what physicists call bosonic fields.
Our purpose here is to show how the notion of chain extension may be used to
reformulate ``local BRST theory" in a setting involving bosonic and
fermionic fields. Physicists and mathematicians often employ different
notation to
describe field theories involving both bosonic and
fermionic fields. We will adopt a hybrid notation which permits us to
deal with the
issues of concern to us.
For  more complete details regarding the notation of this section we
refer the
reader to
\cite{BFLS}. A more systematic development of the properties of
infinite jet
bundles may be found in \cite{KV}.

Let $M$ be a $n$-dimensional manifold and let $\pi: E\longrightarrow{M}$
be a vector
bundle of fibre dimension
$k$ over $M$. Let $J^{\infty}E$  be the infinite jet bundle over $\pi$ with
$\pi_M^{\infty}: J^{\infty}E\longrightarrow{M}$. The bosonic fields under
consideration will be identified with sections of $\pi.$

Fermionic fields are usually modeled by fields called ghost fields. In
order to
describe
the ghost fields we define a finite-dimensional Grassmann algebra ${\cal
G}$ with $n$
linearly independent generators $\{e_a\}.$  Thus a basis of ${\cal G}$
is the set of
products
$e_{a_1}e_{a_2}\cdots e_{a_k}$ where $1\leq a_1 < a_2 < \cdots <a_k \leq
n.$ The
linear
space generated by $e_{a_1}e_{a_2}\cdots e_{k}$ for fixed $k$ is denoted
${\cal G}_k$
and we say that elements of ${\cal G}_k$ have ghost number $k.$ A single
ghost field
is a mapping from $M$ into ${\cal G}_1.$ The Lagrangian of the ``fermionic
sector"  of a Lagrangian field theory is a function of a finite number
$N$ of such
ghost
fields and their derivatives. More precisely such a Lagrangian  is a
function of mappings from $M$ into
${\cal G}_1^N$ and derivatives of such mappings. The values of these
Lagrangians
are in ${\cal G}$ itself as they are usually sums of products of the
fields and
their derivatives.

Similarly, to describe anti-ghost fields we choose a
symmetric algebra ${\cal A}$  which is generated by $n$ linearly
independent elements $\{f_b\}$ of ${\cal A}.$ As before a basis of
${\cal A}$ is the
set of all products  $f_{b_1}f_{b_2}\cdots f_{b_k}$ where
$1\leq b_1 < b_2 < \cdots <b_k \leq n.$ As before
${\cal A}_k$ denotes the subspace generated by products
$f_{b_1}f_{b_2}\cdots f_{b_k}$
and we say that elements of this subspace have anti-ghost number $k.$

Define a vector bundle
by setting $\tilde E_x= E_x\oplus {\cal G}^N_1 \oplus {\cal A}^N_1$ for
each $x\in
M.$ Then $\tilde E$ is a vector bundle over $M$  and the fields under
consideration
will be sections of this bundle.

For simplicity we will assume that the fiber bundle $E\longrightarrow M$
is trivial
and has fiber the vector space $V.$ It follows that $\tilde \pi:\tilde
E\longrightarrow M$ is also trivial with fiber $V\oplus  {\cal G}^N_1
\oplus {\cal
A}^N_1.$ It follows that sections $\psi$ of  $\tilde \pi$ may be
identified with
maps from $M$ into the fiber $V\oplus  {\cal G}^N_1 \oplus {\cal
A}^N_1.$ Moreover
we choose a basis of $V$ once for all so that mappings from $M$ into the
fiber
$V\oplus  {\cal G}^N_1 \oplus {\cal A}^N_1$ are uniquely determined by their
components. We write $\psi(x)=(\phi^i(x),\eta^a(x),P_b(x))$ to denote a
typical
section of $\tilde \pi.$
Observe that the vector components of $\eta$ anti-commute,
$$\eta^a\eta^b=-\eta^b\eta^a,$$ but their scalar
components $\eta^{a \alpha},\eta^{b \gamma}$ are in $C^{\infty}(M)$ and
thus commute.
If $(x^{\mu})$ are the components of a local chart on $M,$ then the
infinite jet of a
section $\psi$  of $\tilde \pi$ may be denoted as
$$j^{\infty}\psi(x)=(x,\phi(x),\eta(x),P(x),\partial_I\phi(x),\partial_J\eta(x),\partial_KP(x)).$$
Here $I,J,K$ are symmetric multi-indices where, for example,
$\partial_I=\partial_{{\mu}_1}\partial_{{\mu}_2}\cdots
\partial_{{\mu}_r}$ are
partial derivatives relative to the chart  $(x^{\mu}).$  We denote
coordinates on
$J^{\infty}\tilde E$ as follows:
$$x^{\mu}(j^{\infty}\psi(p))=x^{\mu}(p), \quad
u_I^i(j^{\infty}\psi(p))=\partial_I\phi^i(x)$$
$$\mu_J^{a\alpha}(j^{\infty}\psi(p))=\partial_J\eta^{a\alpha}(x), \quad
\nu_K^{b\beta}(j^{\infty}\psi(p))=\partial_KP_{b\beta}(x).$$

One may now consider the bivariational-complex of the infinite jet bundle
$J^{\infty}\tilde E.$
The de Rham complex of differential forms $\Omega^*(J^{\infty}\tilde
E,d)$ on
$J^{\infty}\tilde E$ possesses a differential ideal, the ideal ${ C}$ of
contact
forms $\theta$ which satisfy $(j^{\infty}\phi)^* \theta=0$ for all
sections $\phi$
with compact support. This ideal is generated by the contact one-forms,
which in
local coordinates assume the form $\theta^i_I=du^i_I-u^i_{jI}dx^j$
,$\theta^{a\alpha}_J=d\mu^{a\alpha}_J-\mu^{a\alpha}_{jJ}dx^j$,$\theta^{b\beta}_K=d\nu^{b\beta}_K-\nu^{b\beta}_{jK}dx^j.$
Contact one-forms of order $0$ satisfy $(j^{1}\phi)^*(\theta)=0$. For
example, for
the {\it bosonic components} of a field, contact forms of {\it order
zero}  assume the form
$\theta^a= du^a-u^a_j dx^i$ in local coordinates (the muti-index $I$ is
empty in this case).

Using the contact forms, we see that the complex
$\Omega^*(J^{\infty}\tilde E,d)$
splits as a bicomplex $\Omega ^{r,s}(J^ \infty \tilde E)$
(though the finite
level complexes $\Omega^*(J^p\tilde E)$ do not), where $\Omega ^{r,s}(J^
\infty
\tilde E)$
denotes the space of differential forms on $J^\infty \tilde E$ with $r$
horizontal
components and $s$ vertical components. The bigrading is described
by writing a differential $p$-form $\alpha=\alpha_{IA}^{\bf
J}(\theta^A_{\bf J}
\wedge dx^I)$ as an element of
$\Omega^{r,s}(J^{\infty}\tilde E)$, with $p=r+s$, and \begin{eqnarray}
dx^I=dx^{i_1}\wedge...\wedge dx^{i_r} \quad {\rm and}
\quad \theta^A_{\bf J}=\theta^{a_1}_{J_1}\wedge...\wedge \theta^{a_s}_{J_s},
\end{eqnarray} where $\theta^{a_j}_{J_j}$ denotes any one of the three
types of
contact forms defined above.
We are interested
in the complex $$ 0 \to \Omega ^{0,0}(J^\infty \tilde E) \to \Omega
^{1,0}(J^\infty
\tilde E) \to \cdots \to \Omega ^{n-1,0}(J^\infty \tilde E) \to
\Omega ^{n,0}(J^\infty \tilde E)$$ with the differential $d_H$ defined
by $d_H dx^i
D_i$ where
$$D_j=\frac{\partial}{\partial{x^j}}+u_{jI}^{i}\frac{\partial}{\partial
u_{I}^{i}}+\mu_{jJ}^{a\alpha}\frac{\partial}{\partial\mu_{J}^{a\alpha}}+
\nu_{jK}^{b\beta}\frac{\partial}{\partial\nu_{K}^{b\beta}}.$$
Here if $\alpha = \alpha _I dx^I$ then $d_H \alpha = D_i \alpha _I (dx^i
\wedge dx^I)$.
Notice that this complex is exact by the algebraic Poincare lemma.
The algebraic Poincare lemma is essentially the Poincare lemma for the
horizontal
complex on the infinite jet bundle. It is discussed and proved in
\cite{Brandt}, for example.

\begin{definition} To say that $F$ is a {\bf local function } (on
$J^{\infty}\tilde E$) means that
$F$ is a function from $J^{\infty}\tilde E$ into ${\bf R}\otimes {\cal
G}\otimes
{\cal A}$ such that
$F$ factors through the projection of $J^{\infty}\tilde E$ onto
$J^p\tilde E$ for
some nonnegative integer p. Moreover for each section $\psi$ of $\tilde
\pi,$
$\psi=(\phi,\eta,P),$  we require that
$$F(j^{\infty}\psi(x))=F^{{\bf J}{\bf K}}_{AB}(j^{\infty}\phi(x))
(\partial_{J_1}\eta^{a_1})(x)\cdots (\partial_{J_r}\eta^{a_r})(x)
(\partial_{K_1}P_{b_1})(x) \cdots \partial_{K_s}P_{b_s})(x).$$
where $A=\{a_1a_2\cdots a_r\},$ $B=\{b_1b_2\cdots b_s\}$ are
multi-indices and
${\bf J}=\{J_1J_2\cdots J_r\},$   ${\bf K}=\{K_1K_2\cdots K_s\}$ are
vectors of
multi-indices. Here the set of real-valued functions $\{F^{{\bf J}{\bf
K}}_{AB}\}$
are smooth functions on the jet bundle $J^{\infty}E.$
We denote this algebra of local functions on $J^{\infty}\tilde E$ by
$Loc=Loc_{\tilde E}.$
\end{definition}

\begin{definition}
A  {\bf local functional}
\begin{eqnarray}
&& {\cal L} [\psi]=\int_M L(x,\phi^{(p)}(x),
,\partial_{j_{1}}\partial_{j_{2}}\cdots \partial_ {j_{s}}
\eta^{b_j}\partial_{i_{1}}\partial_{i_{2}}\cdots
\partial_{i_{r}}P^{a_i})dvol_M
\end{eqnarray}
is the integral over $M$ of a local function $L$ evaluated on sections
$\psi$ of
$\tilde E$ of compact support. Thus
$${\cal L}(\psi)=\int_M j^{\infty}(\psi)^*(L) dvol_M$$
for some local function $L\in Loc.$
\end{definition}

{\bf Remark.} Notice that the integral in the last definition is the
integral of a
vector-valued function and the result is then a vector in the space
${\bf R}\otimes
{\cal G}\otimes {\cal P}.$ We presume that the integral respects the
grading on this
algebra as this is usually the case in quantum field theory in which
integrals of
this type regularly arise.
\bigskip

In physical applications, the BRST operator is often a mapping  $S:
Loc\longrightarrow Loc$ which commutes with
the horizontal differential $d_H$ and which satisfies the condition
$S^2=d_H \tilde k=\partial_{i }k^{i }$ where $\tilde
k=k^i(-1)^{i-1}[\frac{\partial}{\partial x^i}\lh (dx^1\wedge dx^2\wedge
\cdots
dx^n)]$ ( see \cite{Brandt}
and \cite{ HT}). Thus
in the  Lagrangian field setting,the BRST  differential is essentially a
differential  on a quotient  space and consequently on the space $\cal F$
of local functionals.
Notice that  $(S,\cal F)$ is a $D$-algebra, and if we set ${\cal
B}=d_H(\Omega^{n-1,0})$= divergences, then we can obtain a local BRST
theory via a
chain
extension  for the triple $(\cal B,\cal F,S)$.To apply our theorem, we
need to
construct a homological resolution for space $\cal F$.

Obviously,  $\cal F\simeq$ $ H^n(\Omega^{*,0},d_{H})$ and $\cal B=$
$d_H(\Omega^{n-1,0}).$ By the  algebraic Poincare lemma
\cite{Brandt}(this is
essentially the Poincare lemma for the horizontal differential), we have a
commutative diagram:
$$\begin{array}{ccccccccccc}\\& & & d_H & & d_H & & d_H & & d_H&\\ {\bf
    R}& \longrightarrow &\Omega^{0,0}&\longrightarrow &\dots &
\longrightarrow &
\Omega^{n-2,0} &\stackrel {\textstyle}{\longrightarrow}
&\Omega^{n-1,0}&\stackrel {\textstyle}{\longrightarrow}
&\Omega^{n,0}\\ \downarrow\eta & &\downarrow\eta & & & &
\downarrow\eta & & \downarrow\eta & & \downarrow\eta\\ 0 &
\longrightarrow & 0 &\longrightarrow &\dots & \longrightarrow
&0&\longrightarrow & 0 &\longrightarrow & H_0={\cal F}
\end{array}$$
Thus $(\Omega^{*,0},d_H)$ provides a homological resolution for $\cal
F$. It follows
that if there exists an appropriate contracting homotopy then we can use
Theorem 3
to construct an operator
$\tilde S=d_H+S+l_3$ which serves in place of the BRST operator $S$ in
the usual
approach to local BRST theory. The explicit expression
for $\tilde S$ is interesting but more complicated because the
contracting homotopy
given by the algebraic Poincare lemma is not so simple (see \cite{Olver}
and \cite
{Brandt}).

\subsection{Consistent deformation}
The purpose of this section is to reanalyze the long-standing problem of
constructing
consistent interactions among fields with a gauge freedom using the
formalism
of chain extensions. Analysis of this problem has a long history which we
will not reiterate here, preferring instead to direct attention to two
papers.
The first of these is the paper by Berhends, Burgers, and Van Dam
\cite{BBvD} where there is a frontal assault on the problem using direct
methods
and the second is the paper by Barnich and Henneaux \cite{BH} where the
problem is addressed using  the antibracket formalism. The latter paper
reformulates the problem as a deformation problem
in the sense of deformation theory \cite{G}, namely that of
{\em{deforming
consistently the master equation}}. They show, by using the properties of
the antibracket,
that there is no obstruction to constructing interactions that consistently
preserve
the gauge symmetries of the free theory if one allows the interactions to be
non local.

It is our intent here show that the notion of chain extension clarifies the
description of the algebraic aspects of the anti-bracket approach  of
the problem.
Recall that  if $L$ is a Lagrangian with gauge degrees of freedom, then the
simultaneous deformation of the Lagrangian and it's gauge transformations
$\delta_\beta$ can be related to the deformation of the master equation.

In the Batalin-Vilkovisky formalism for gauge theories, which we
consider to be irreducible, one introduces,
besides the original fields $\phi^i$ of ghost number $0$
and the ghosts $C^\alpha$ of ghost number $1$ (related
to the gauge invariance), the corresponding antifields $\phi^*_i$ and
$C^*_\alpha$ of opposite Grassmann parity and ghost number $-1$ and
$-2$ respectively \cite{HT}.
An  antibracket is defined by
requiring that the fields $\phi^A\equiv(\phi^i,C^\alpha)$ and
antifields $\phi^*_A$ be conjugate~:
\begin{eqnarray}
(\phi^A(x),\phi^*_B(y))=\delta^A_B\delta(x-y)
\end{eqnarray}
Physicists (see \cite{HT} ) then define the anti-bracket  for arbitrary
functionals
$A_1$ and $A_2$
of these  ``extended fields" by
\begin{eqnarray}
(A_1,A_2)=\int ({\delta^RA_1\over\delta\phi^A}
{\delta^LA_2\over\delta\phi^*_A}
-{\delta^RA_1\over\delta\phi^*_A}
{\delta^LA_2\over\delta\phi^A})d^nx
\end{eqnarray}
where the operator $\frac  {\delta^R}{\delta \phi^A}$ is a functional
derivative
acting from the right and the other terms are defined similarly.

Notice that functionals such as $A_1,A_2$ above are not functions of the
fields
and antifields but rather are integrals of such functions. A more
precise formulation
of these concepts was provided in the last subsection on local functionals.

In  the following paragraphs,
the parity $\epsilon_A$ of a functional $A$ is defined by requiring that
$\epsilon_A=1$
when the ghost number of $A$ is odd and $\epsilon_A=0$ when the ghost
number of $A$ is even.
The central goal of the formalism is the construction of a proper
solution to the master equation
\begin{eqnarray}
(S,S)=0,
\end{eqnarray} where $\epsilon_S=0.$
The functional $S$ is required to begin with the classical action
$S_0,$  to
which one couples
via antifields the gauge transformations with the gauge parameters
replaced by the ghosts. The BRST symmetry is then canonically generated
by the
antibracket
through the equation~:
\begin{eqnarray}
s=(S,\cdot).
\end{eqnarray}
The description of the anti-field formalism we have given here closely
follows that
of \cite{HT}.

Thus we require that the components of both the fields
$\phi^A$ and anti-fields $\phi^*_A$ be mappings from a given manifold
$M$ into the
complex field.  We then define
${\cal R}$ to be the algebra of all  functionals $A=A(\phi^B,\phi^*_B)$
of the fields,
anti-fields, and their derivatives.  If $A_1,A_2$ are
in  ${\cal R}$ then the anti-bracket $(A_1,A_2)$ defined above is then
again an element
of
${\cal R}.$ The bracket defined in this manner is odd, carries ghost
number +1, and
obeys
the (graded) Jacobi identity. More explicitly, if $\epsilon_F$ denotes
the parity of
an element $F$ of ${\cal R}$ we have for  $A,B,C \in {\cal R}:$
$$(A,B)=-(-1)^{(\epsilon_A+1)(\epsilon_B+1)}(B,A)$$
$$(-1)^{(\epsilon_A+1)(\epsilon_C+1)}(A,(B,C))+(-1)^{(\epsilon_B+1)(\epsilon_A+1)}(B,(C,A))
\quad \quad \quad \quad$$
$$\quad \quad \quad \quad \quad \quad \quad
\quad+(-1)^{(\epsilon_C+1)(\epsilon_B+1)}(C,(A,B))=0.$$

For our exposition below we find it useful to extend the definition of the
anti-bracket
to the space of formal power series  ${\cal
R}[[t]]=\{\sum_{i=0}^{\infty}a_it^i | \quad
a_i\in {\cal R} \}.$
We define an anti-bracket on ${\cal R}[[t]]$ by:
$$(\sum_{i=0}^{\infty}a_it^i,\sum_{j=0}^{\infty}b_jt^j)=\sum_{i,j=0}^{\infty}(a_i,b_j)t^{i+j}.$$
Clearly, ${\cal R}[[t]]$ is both a vector space over the field $k$ of
complex numbers
and a module over the algebra $k[[t]].$

\bigskip

{\em{It will be fairly obvious from the exposition below that in our
formalism
${\cal R}$ could be any vector space over
an arbitrary  field $k$ of characterisic  zero provided one has an
anti-bracket
defined on
${\cal R}$ subject to the conditions above.}}
\bigskip

We now proceed to analyze the master equation
using chain extensions.
Let $S_0\in {\cal R}$ be a solution of the master equation
\begin{eqnarray}
(S_0,S_0)=0
\end{eqnarray}
and let $s_0=(S_0,\cdot )$ denote the corresponding BRST differential.
The problem of deforming  $S_0$ is the problem of finding a formal power
series
$S(t)=\sum_0^{\infty}S_i t^i$ in ${\cal R}[[t]]$ such that
$(S(t),S(t))=0$ and such
that
$\epsilon_{S_i}=0$ for all $i.$  Expansion of
the master equation in the deformation parameter gives:
\begin{eqnarray}
(S_0,S_0)=0\nonumber\\(S_0,S_1)=0\nonumber\\(S_1,S_1)+(S_0,S_2)+(S_2,S_0)=0\nonumber\\
{\cdots} \nonumber
\end{eqnarray}
In general, $\sum_{i+j=n}(S_i,S_j)=0$.
In order to solve for the $S_n$ inductively, we notice that:
\begin{eqnarray}
0=\sum_{i+j=n,i,j\geq 0}(S_i, S_j)=\sum_{i+j=n,i,j\geq 1}(S_i,
S_j)+(S_0,S_n)+(S_n,S_0)
\end{eqnarray}
Since the anti-bracket is graded commutative and gh$(S)=0$, we have
$(S_0,S_n)=(S_n,S_0)$,thus:
\begin{eqnarray}
\sum_{i+j=n,i,j\geq 1}(S_{i}, S_{j})=-2(S_{0},S_{n})
\end{eqnarray}
Set $ R_n=\sum_{i+j=n,i,j\geq 1}(S_i ,S_j),$ obviously $R_n$ is
determined by the
first $n-1$ terms:\newline $S_1,S_2\cdots S_{n-1}$ of the sequence $S$
and $R_n$ is a  $s_0$-coboundary, which we call the   $n$-th obstruction
to the
deformation
of $S_0.$  By definition the $n$-th order deformation of $S_0$ is a sequence
$S_1,S_2\cdots , S_{n}$ in ${\cal R}$  such that
\begin{eqnarray}
(\sum_0^{n}S_i t^i,\sum_0^{n}S_j t^j)\equiv 0 \quad (mod \quad  t^{n+1}).
\label{S} \end{eqnarray}

Given an $n$-th order deformation of $S_0$ clearly the obstruction to
obtaining an
$(n+1)$-th order deformation of $S_0$  is $R_{n+1}=\sum_{i+j=n+1,i,j\geq
1}(S_i,
S_j).$

Define an submodule $I$ in ${\cal R}[[t]]$ by
$I={\cal R}[[t]] t^{n+1}=\{\sum_{i\geq n+1}a_i t^i | \quad a_i\in {\cal
R}\};$ thus we
have a short exact sequence:

\begin{eqnarray}
& &  0\longrightarrow I\longrightarrow {\cal R}[[t]]\longrightarrow {\cal
R}[[t]]/(t^{n+1})\longrightarrow 0.
\end{eqnarray}

Define a $k$-linear mapping $D$ on ${\cal R}$ with values in  ${\cal
R}[[t]]$ by
$D(x)=(\sum_0^{n}S_i t^i,x)$ for each $ x\in {\cal R}.$ The scalar
multiplication of the field $k$ on ${\cal R}[[t]]$ can clearly be
extended to a module multiplication of the power series ring $k[[t]]$ on
${\cal R}.$
Moreover the operator $D$ can  be extended to  a $k[[t]]$ module
homomorphism
of  ${\cal R}[[t]]$ with values in ${\cal R}[[t]].$
Notice that $D(I)\subseteq I$ and as a consequence there is an induced
operator on
$\bar {\cal R}[[t]]= {\cal R}[[t]]/I$
which we denote by$\bar D.$ We will show below that $D^2(I)\subset I$ and
consequently that $\bar D^2=0.$ Thus we obtain a triple $(I,\bar {\cal
R}[[t]], \bar
D);$ we intend to construct a chain extension for this
triple.

To do this we first  introduce a ``superpartner set of ${\cal R}$" which
we  denote
by ${\cal R}[1]$ and which is defined by ${\cal R}[1]=\{a^* | \quad
a\in {\cal R}\}$
where
$a^*\leftrightarrow a$ is a one to one correspondence and
$\epsilon(a^*)=\epsilon(a)+1.$

Let $X_0={\cal R}[[t]]$ and $X_1={\cal R}[1][[t]]t^{n+1}=\{\sum_{i\geq
n+1}a_i^* t^i|
\quad a_i^*\in {\cal R}\},$ so that one has  a homological resolution
of  ${\cal R}[[t]]/I:$
\begin{eqnarray}
& & 0\longrightarrow X_1 \stackrel{l_1}{\longrightarrow} X_0 \longrightarrow
0.\label{z}
\end{eqnarray}
Here $l_1:X_1$ $\longrightarrow X_0 $ is defined by
$l_1(x)=\sum_{i\geq n+1}a_it^i$ for each $x$ such that $ x=\sum_{i\geq
n+1}a_i^*t^i.$
We can now construct a contracting homotopy similar to earlier
definitions. In order
to avoid confusion, we denote the contracting homotopy in that which
follows by $h.$
Since the sequence (\ref{z}) is split, there exists a natural direct sum
decomposition
of $X_0=I\oplus X_0/I.$ Define ${\cal F}=X_0/I,$  $\eta:
X_0\longrightarrow \cal F$
by $\eta =proj|_{\cal F}$ and $ \lambda $ by $\lambda =i_{{\cal F}
\rightarrow X_0}$
(the inclusion mapping).
Finally, since $X_0=I\oplus \cal F,$  define the contracting homotopy
$h$ by
$h(x)=-\sum_{i\geq n+1}a_i^*t^i$ for each $x\in I,x=\sum_{i\geq n+1}a_i
t^i,$ and
define $h$ to be
zero on ${\cal F}.$
It's easy to show that
$\lambda \circ \eta -1_{X_*}=l_1\circ h +h \circ l_1,$ thus we may
define a chain
extension for the triple $(I,{\cal F},\bar D)$ via the proposition below.

Let $S_D=\sum_{i=0}^nS_it^i\in X_0={\cal R}[[t]]$ so that $D:X_0
\rightarrow X_0$ is
given by $D(x)=(S_D,x) $ for all $x\in X_0.$ It is easy to show that
$$(S_D,(S_D,x) )= \frac {1}{2}((S_D,S_D),x)$$
for all $x\in X_0$ and that $D^2(I)\subseteq I.$
It follows that one has an induced operator $\bar D$ on ${\cal
R}[[t]]/I.$ Using this
fact we obtain the following proposition.

\begin{proposition}  Assume that  $ {\cal R}$ is a graded linear space
(graded with respect
to ghost number above)
and that it has an anti-bracket structure which is graded commutative
and which
satisfies the graded Jacobi
identity as above. Let ${\cal R}[[t]]$ denote the corresponding space of
formal power series
and let $I$ denote the $k[[t]]$ submodule ${\cal R}[[t]]t^{n+1}.$  If
$S_0,S_1, \cdots
,S_n$ are even elements of ${\cal R}$ and $D$ is the operator
on ${\cal R}[[t]]$ defined by $D(x)=(\sum_{i=0}^nS_it^i,x), \quad x\in
{\cal R}[[t]],$
then $D(I)\subseteq I $ and the induced operator $\bar D$ on ${\cal
R}[[t]]/I$ is
nilpotent. Moreover the triple $(I,{\cal R}[[t]]/I,\bar D)$ admits a
chain extension with chain maps $l_1,l_2,l_3$ such that the map
$S=l_1+l_2+l_3$
satisfies $S^2=0.$ The
explicit homological resolution $X_*$ is defined in the proof below.
\end{proposition}

\begin{proof} Recall that $X_0={\cal R}[[t]], X_1={\cal R}[1][[t]]t^{n+1},$
and define a
mapping
$ l_2:X_0 \rightarrow X_0$ by $ l_2(x)=D(x)=(\sum_{i=0}^nS_it^i,x), x\in
X_0.$
Since we intend to apply Theorem 3, recall that in that context, ${\cal
B}=l_1(X_1)$
and consequently, ${\cal B}=I$ in the present case. Thus $l_2({\cal B})=
l_2({\cal R}[[t]]t^{n+1})=D({\cal R}[[t]]t^{n+1}) \subseteq {\cal
R}[[t]]t^{n+1}=I={\cal B}.$ This is one of the conditions required to
apply Theorem 3.
Also note that for $x\in X_0,$
$ l_2^2(x)=D^2(x)=\frac {1}{2}((S_D,S_D),x)$ and by (\ref{S})
$((S_D,S_D)=0 \quad (mod(t^{n+1}).$
Thus $ l_2^2(X_0)\subseteq I={\cal B}.$

Thus the hypothesis of Theorem 3 holds and consequently the
theorem asserts
the existence of maps
$ l_1, l_2, l_3$ which define a chain extension, in particular their sum
has square
zero.
\end{proof}

From the proposition we know that the triple $(I,{\cal R}[[t]]/I,\bar
D)$ admits a
chain extension with chain maps $l_1,l_2,l_3$ defined by Theorem 3. In
applications, it
is useful to have explicit formulas for these maps. Each mapping $l_i$
is a mapping
from
$X_*$ to $X_*$ which is uniquely determined by its values on ${\cal R}$
and ${\cal
R}[1].$ Our next theorem provides explicit formulas which show how to
obtain the maps
$l_i,i=1,2,3.$

\begin{theorem}
The chain maps $l_1,l_2,l_3,$ which define the chain extension $S=
l_1+l_2+l_3,$
guaranteed by the proposition above, are uniquely determined by their
values on the
linear spaces  ${\cal R}$ and
${\cal R}[1]$ and are given in a concise form as follows. The mapping
$l_1$ is
determined by
its values on $ {\cal R}[1]$ with $ l_1(a^*)=a.$ The map $l_2$ is given
by its
values on
${\cal R}$ and ${\cal R}[1],$ respectively by
$$ l_2(a)=\sum_{j=0}^{n}(S_i,a)t^i, \quad  l_2(a^*)= -
\sum_{j=0}^{n}(S_i,a)^*t^i.$$
Finally, $l_3$ is uniquely determined by its values on ${\cal R},$
$$ l_3(a)=-\frac{1}{2}\sum_{n+1\leq i+j\leq 2n}((S_i,S_j),a)^* t^{i+j}.  $$

\end{theorem}
\begin{proof} By the proposition above and Theorem 3 we have,
$$
l_2(a^*)=hl_2l_1(a^*)=hl_2(a) =h(\sum_0^{n}S_i t^i,a)=h(\sum_0^{n}(S_i,a)
t^i)
=-\sum_0^{n}(S_i,a)^* t^i.$$

Before evaluating $l_3$ ,we replace $l_2$ by $D,$ which appears in the
canonical
transformation and in the corollary above.
Then by Theorem 3, $l_3(a)=hl_2^2(a)=h(D,(D,a))$.
Since $$ D(x)=(\sum_{i=0}^{n}S_i t^i,x),$$ we have:

\begin{eqnarray}
l_3(a)=\frac{1}{2}h((\sum_{i=0}^{n}S_i t^i,\sum_{j=0}^{n}S_j
t^j),a)\nonumber\\
=\frac{1}{2}h(\sum_{k=0}^{2n}\sum_{i+j=k}((S_i,S_j),a)t^{i+j})\nonumber\\
=-\frac{1}{2}\sum_{k=n+1}^{2n}\sum_{i+j=k}((S_i,S_j),a)^*
t^{i+j})\nonumber\\
=-\frac{1}{2}\sum_{n+1\leq i+j\leq 2n}((S_i,S_j),a)^* t^{i+j} \nonumber \\
=-\frac{1}{2}R_{n+1}
-\frac{1}{2}\sum_{n+2\leq i+j\leq 2n}((S_i,S_j),a)^* t^{i+j} .
\end{eqnarray}
This finishes the proof of the theorem.
\end{proof}
\bigskip

{\bf Remark.} Observe that the last calculation above identifies the
obstruction $R_{n+1}$ as a summand of the value of $l_3$ on ${\cal R}.$
In particular, if $l_3$ is zero then the obstruction
to further deformation vanishes.
\bigskip

If one wishes to consider infinitesimal deformations, simply take $n=1,$
then by the
theorem, we have:
$$S(a)=S_0(a)+S_1(a)t-((S_1,S_1),a)t^2$$
$$S(a^*)=a-(S_0,a)^*-(S_1,a)^* t.$$
Extend these as derivations and we obtain a consistent deformation of
the space of fields $\cal R $
by enlarging $\cal R $ to the space ${\cal R}[1] $,or physically
speaking, the deformation
is achieved  by adjoining the superpartners to the space of original fields.

\subsection{Deformation theory, sh-Lie algebras }

In the last two decades,  deformations of various types of structures
have assumed an ever
increasing role in mathematics and physics. For each such deformation
problem  a goal  is to determine
if all related deformation obstructions  vanish and many beautiful
techniques been
developed  to determine when this is so. Sometimes genuine deformation
obstructions arise and occasionally  that closes mathematical
development in such cases, but
in physics such problems are dealt with by
introducing  new auxiliary fields to kill such obstructions. This idea
suggests that one might deal with deformation problems  by enlarging the
relevant category to a new category obtained by
appending additional algebraic structures to the old category.

In the present subsection, we consider deformations of Lie algebras.
In order to be complete we review  basic facts on Lie
algebra deformations; more detail may be found in the book edited by
M.Hazawinkel
and M.Gerstenhaber \cite{HG}.

Let $A$ be a $k$-algebra and $\alpha$ be its multiplication, i.e.,
$\alpha$ is  a
$k$-bilinear map $A\times A \longrightarrow A$ defined by
$\alpha(a,b)=ab.$ A deformation of $A$ may be defined to be a formal
power series
$\alpha_t=\alpha+t\alpha_1+t^2\alpha_2+\cdots $
where each $\alpha_i:A\times A \longrightarrow A$ is a $k$-bilinear map
and the
``multiplication" $\alpha_t$ is formally of
the same "kind" as $\alpha,$ e.g., it is  associative or Lie or whatever
is required. One
technique used  to set up a deformation problem is to extend a
k-bilinear mapping
$\alpha_t:A\times A \longrightarrow A[[t]]$ to  a  $k[[t]]$-bilinear
mapping $\alpha_t:A[[t]]\times A[[t]]
\longrightarrow A[[t]].$ A mapping
$\alpha_t:A[[t]]\times A[[t]] \longrightarrow A[[t]]$ obtained in this
manner  is
necessarily uniquely determined by it's
values on $A\times A.$ In fact we would not regard the mapping
$\alpha_t:A[[t]]\times A[[t]] \longrightarrow A[[t]]$
to be a deformation of $A$ unless it is determined by it's values on $
A\times A$.

From this point on, we assume that $(A,\alpha)$ is a Lie algebra,i.e.,
we assume that
$\alpha(\alpha(a,b),c)+\alpha(\alpha(b,c),a)+\alpha(\alpha(c,a),b)=0$.
Thus the problem of
deforming a Lie algebra $A$ is equivalent to the problem of finding a
mapping
$\alpha_t:A\times A \longrightarrow A[[t]]$ such that
$\alpha_t(\alpha_t(a,b),c)+\alpha_t(\alpha_t(b,c),a)+\alpha_t(\alpha_t(c,a),b)=0.$
If we set $\alpha_0=\alpha$ and expand this Jacobi identity by making
the substitution
$\alpha_t=\alpha+t\alpha_1+t^2\alpha_2+\cdots ,$ we get the equation
\begin{eqnarray}
\sum_{i,j=0}^{\infty}[\alpha_j(\alpha_i(a,b),c)+\alpha_j(\alpha_i
(b,c),a)+\alpha_t(\alpha_i (c,a),b)]t^{i+j}=0
\end{eqnarray}
and consequently a sequence of deformation equations;
\begin{eqnarray}
\sum_{i,j\geq 0,i+j=n}[\alpha_j(\alpha_i(a,b),c)+\alpha_j(\alpha_i
(b,c),a)+\alpha_t(\alpha_i (c,a),b)]=0.
\end{eqnarray}
The first two equations are:
\begin{eqnarray}
\alpha_0(\alpha_0(a,b),c)+\alpha_0(\alpha_0(b,c),a)+\alpha_0(\alpha_0(c,a),b)=0
\\
\alpha_0(\alpha_1(a,b),c)+\alpha_0(\alpha_1(b,c),a)+\alpha_0(\alpha_1(c,a),b)+\alpha_1(\alpha_0(a,b),c)\nonumber
\\
  +\alpha_1(\alpha_0(b,c),a)+\alpha_1(\alpha_0(c,a),b)=0
\end{eqnarray}

We can reformulate the discussion above in a slightly more compact form.
Given
a sequence
$\alpha_n: A\times A\longrightarrow A$
of bilinear maps, we define ``compositions" of various of the
$\alpha_n$ as follows:
\begin{eqnarray}
& & \alpha_i\alpha_j:A\times A \times A \longrightarrow A
\end{eqnarray}
is defined by
\begin{eqnarray}
& & (\alpha_i \alpha_j)(x_1,x_2,x_3)=\sum_{\sigma \in
unsh(2,1)}(-1)^{\sigma} \alpha_i (\alpha_j(x_{\sigma(1)},
x_{\sigma(2)}),x_{\sigma(3)})
\end{eqnarray}
for arbitrary $x_1,x_2,x_3\in A.$

Thus the deformation equations are equivalent to following equations:
\begin{eqnarray}
\alpha_{0}^{2}=0 \\ \alpha_0\alpha_1+\alpha_1\alpha_0=0 \\
\alpha_1^2+\alpha_0\alpha_2+\alpha_2\alpha_0=0\\
\cdots \nonumber\\ \Sigma_{i+j=n} \alpha_i\alpha_j=0\\ \cdots.
\end{eqnarray}

Define a bracket on the sequence $\{\alpha_n\}$ of mappings by
$[\alpha_i,\alpha_j]=\alpha_i\alpha_j+\alpha_j\alpha_i$ and a
``differential"
$d$ by $d=ad_{\alpha_0}=[\alpha_0, \cdot ],$ the ``adjoint
representation" relative
to $\alpha_0.$  Notice that the second equation in the list above is
equivalent to the
statement that $\alpha_1$ defines a cocycle $\alpha_1\in Z^2(A,A)$ in
the Lie
algebra cohomology of $A.$
Moreover it is known that the second cohomology group $H^2(A,A)$
classifies the equivalence class of infinitesimal deformations of $A$
\cite{HG}.
This being the case we refer to the triple $(A,\alpha_0,\alpha_1)$ as
being initial
conditions for deforming the Lie algebra $A.$
Notice that the third equation in the above list can be rewritten as
\begin{eqnarray}
[\alpha_1,\alpha_1]=-[\alpha_0,\alpha_2]=-d\alpha_2
\end{eqnarray}
When this  equation holds one has then that  $[\alpha_1,\alpha_1]$ is a
coboundary
and so defines
the trivial element of $H^3(A,A)$ for any given deformation $\alpha_t.$
Thus if
$[\alpha_1,\alpha_1]$ is not a coboundary,  then we may regard
$[\alpha_1,\alpha_1]$ as the first  obstruction to deformation and in
this case we
can not deform
$A$ at second order. In general, to say that there exists a  deformation of
$(A,\alpha_0,\alpha_1)$ up to order $n-1,$ means that  there
exists a sequence of maps $\alpha_0,\cdots, \alpha_{n-1}$ such that
$\sum_{\sigma \in unsh(2,1)}(-1)^{\sigma} \alpha_t  (\alpha_t(x_{\sigma(1)},
x_{\sigma(2)}),x_{\sigma(3)})=0  \quad (mod \quad t^n). $
If this is the case and if there is an obstruction to  deformation  at
$n$th order,
then it follows  that $\rho_n=-(\Sigma_{i+j=n,i,j>0}
\alpha_i\alpha_j)$ is in some sense the obstruction and $[\rho_n]$ is a
nontrivial
element of $H^3(A,A).$  If $[\rho_n]\neq 0$,then the process of obtaining a
deformation will terminate at order $n-1$ due to the existence of the
obstruction
$\rho_n.$  In principal, it is possible that one could return to the
beginning and
select different
terms for the $\alpha_i$ but when this fails what can one say? This is
the issue in
the remainder
of this section.

Indeed  the central point of this section  is to show
that when there is an obstruction to the deformation of a Lie algebra,
one can
use
the obstruction itself to define one of the structure mappings of an
sh-Lie algebra.
Without loss of generality, we  consider a  deformation problem which
has a first
order obstruction.

The required sh-Lie structure lives on a graded vector space $X_*$ which
we define
below. This space in degree zero is given by $X_0=A[[t]]=(\{\Sigma
a_it^i |\quad a_i\in A\}.$ The spaces  $\cal B=$
$ <t^2>=$ $A[[t]]\cdot t^2=$ $\{\Sigma_{i\geq 2} a_it^i |a_i\in A\}$ and
$\cal F= $
$X_0/{\cal B}=$  are also relevant to our construction. Notice that
$\cal F$ is isomorphic to $\{a_0+a_1 t | \quad a_0,a_1 \in A\}$ as a
linear space and that   $X_0,\cal B$ are both $k[[t]]$-modules while
$\cal F$ is a $k[[t]]/<t^2>$ module (recall that $k$ is
underlying field of $A$).
To summarize, we have
following short exact sequence:
$$
0{\longrightarrow} {\cal B} {\longrightarrow} {X}_{0} {\longrightarrow
}{\cal F}
\longrightarrow{0}.
$$
Suppose that the initial Lie structure of $A$ is given  by
$\alpha_0: A\times A \longrightarrow A$ and denote a fixed
infinitesimal deformation  by $[\alpha_{1}] \in H^{2}(A,A).$  One of the
structure mappings
of our sh-Lie structure will be determined by the mapping $\tilde l_2 :
X_0 \times X_{0} \longrightarrow X_{0} $
defined as follows: for any $a,b \in A$, let
\begin{eqnarray}
\tilde l_2(a,b)=\alpha_0(a,b)+\alpha_1(a,b)t
\end{eqnarray}
and extend it to $X_0$ by requiring that it be  $k[[t]]$-bilinear.
Obviously, $\tilde
l_2$ induces a Lie bracket $[,]$ on $\cal F$, but if the  obstruction
$[\alpha_1,\alpha_1] $ is not zero,
then $\tilde l_2^2 \neq 0$ and consequently $\tilde l_2$ can not be a
Lie bracket on $X_0$
(since it doesn't satisfy the Jacobi identity).

To deal with this obstruction we will show that we can use
$\alpha_0,\alpha_1$ to
construct an
sh-Lie structure with at most three nontrivial structure maps
$l_1,l_2,l_3$ such
that the value of $l_3$ on $A\times A\times A$
is the same as that of $[\alpha_1,\alpha_1].$ In
particular, $l_3$ will vanish if and only if the obstruction
$[\alpha_1,\alpha_1]$
vanishes. Thus the sh-Lie algebra encodes the obstruction to deformation
of the Lie
algebra $(A,\alpha_0).$

The required sh-Lie algebra lives on a certain homological resolution
$(X_*,l_1)$ of
${\cal F},$
so our first task is to construct this resolution space for $ {\cal F}.$ To
do this let's introduce a ``superpartner set of $A,$" denoted by  $A[1],$ as
follows:
for each $a\in A,$ introduce $a^*$ such that $a^*\leftrightarrow a$ is a
one to one
correspondence and define $\epsilon(a^*)=\epsilon(a)+1.$
Let $X_1=A[1][[t]]t^2$ and define a map $l_1: X_1 \longrightarrow X_0$ by
$$ l_1(x)=\Sigma_{i\geq 2} a_i t^i \in X_0, \quad \quad
x=\Sigma_{i\geq 2}a_i^*t^i \in X_1 .      $$

Notice that this is just the $k[[t]]$ extension of the
$a^*\leftrightarrow a$  map.
Since $l_1$ is injective, we obtain a homological resolution
$X_*=X_0\oplus X_1$
due to the fact that  the  complex defined by:
\begin{eqnarray}
0\longrightarrow X_1\stackrel{l_1}{\longrightarrow}X_0\longrightarrow 0
\end{eqnarray}
has the obvious property that $H(X_*)=H_0(X_*)\simeq \cal F.$

The $sh$-Lie algebra being constructed will have the property that
$l_n =0,n\geq 4.$ Generally sh-Lie algebras can have any number of
nontrivial
structure maps. The fact that all the structure mappings of our sh-Lie
algebra are
zero with the exception of $l_1,l_2,l_3$ is an immediate consequece of
the fact
that we are able to produce a
resolution of the space $\cal F$ such that $ X_k =0 $ for $k \geq 2.$ In
general
such resolutions do not exist and so one does not have $ l_n =0 $ for
$n\geq 4.$

In order to finish the
preliminaries, we now  construct a contracting homotopy $s$
such that following commutative diagram holds:
$$
\begin{array}{ccccc}
&&s&& \\
0\longrightarrow & X_1 &\stackrel
{\textstyle\longleftarrow}{\longrightarrow} &X_0&\longrightarrow 0\\
&&l_1&& \\ 
& \lambda\Bigl\uparrow\Bigr\downarrow\eta & &
\lambda\Bigl\uparrow\Bigr \downarrow\eta & \\
&&&& \\
0\longrightarrow &0&\longrightarrow & \cal F & \longrightarrow 0
\end{array}
$$
Clearly the  linear space $X_0$ is the direct sum of ${\cal B}$ and
a complementary subspace which is isomorphic to ${\cal F};$ consequently
we have
$X_0 \simeq {\cal B} \oplus \cal F.$  Define $\eta=proj\mid_{\cal F},
\lambda=i_{{\cal
F}\rightarrow X_0}$
and a contracting homotopy $s: X_0\longrightarrow X_1$ as follows:
write $X_0={\cal B}\oplus \cal F,$  set $s|_{\cal F}=0,$ and let $
s(x)=-x^*$ for
all $ x\in {\cal B}.$
It is easy to show that
$ \lambda \circ \eta  -1_{X_*}=l_1\circ s +s \circ l_1.$
In order to obtain the sh-Lie algebra referred to above, we apply a
theorem of
\cite{BFLS}. The hypothesis of this theorem requires the existence of a
bilinear mapping
$\tilde l_2$ from $X_0\times X_0$ to $X_0$ with the properties that  for
$c,c_1,c_2,c_3\in X_0$ and $b\in {\cal B}$
$(i) \quad  \tilde l_2(c,b)\in {\cal B}$ and $(ii) \quad  \tilde
l_2^2(c_1,c_2,c_3)\in {\cal B}.$
To see that $(i)$ holds notice that if $p(t),q(t)\in X_0=A[[t]],$ then
$\tilde l_2(p(t),q(t)t^2)=r(t)t^2$ for some $r(t)\in A[[t]]=X_0.$
Also note that the fact that $\tilde l_2$ induces a Lie bracket on ${\cal
F}=X_0/{\cal B}$ implies that
$\tilde l_2^2$ is zero modulo ${\cal B}$ and $(ii)$ follows. Thus $X_*$
supports an
sh-Lie structure
with only three nonzero structure maps $l_1,l_2,l_3$ (see the remark at
the end of
\cite{BFLS}).

\begin{theorem} Given a Lie algebra $A$ with Lie bracket $\alpha_0$ and an
infinitesimal obstruction $[\alpha_1]\in H^2(A,A)$ to deforming
$(A,\alpha_0),$ there
is an sh-Lie algebra on the graded space
$(X_*,l_1)$ with structure maps $\{l_i\}$ such that $l_n=0$ for $n\geq
4.$ The
graded space
$X_*$ has at most two nonzero terms $X_0=A[[t]],X_1=A[1][[t]]t^2.$
Finally, the maps
$l_1,l_2,l_3$ may be given explicitly in terms of the maps
$\alpha_0,\alpha_1.$
\end{theorem}

${\bf Remark:}$ The mapping $l_1$ is simply the differential of the
graded space $(X_*,l_1).$ The mapping $l_2$
restricted
to $X_0 \times X_0$ is the mapping $\tilde l_2$
defined directly in terms of $\alpha_0,\alpha_1$ above.  On $X_1\times
X_0,$ $l_2$ is
determined by
$l_2(a^*t^2,b) =t^2(\alpha_0(a,b)^*+\alpha_1(a,b)^*t)$ for $a^*\in
A[1],b\in A.$
Finally, $l_3$ is uniquely determined by its values on $A\times A\times
A\subset
X_0\times X_0\times X_0$ and is explicitly a multiple of the obstruction
to the
deformation
of $(A,\alpha_0), $ in particular,
$  l_3(a_1,a_2,a_3)=-t^2[\alpha_1,\alpha_1](a_1,a_2,a_3),  a_i\in A.$

\proof { The sh-Lie structure maps are given by Theorem 7 of
\cite{BFLS}. The fact that
$l_n=0,n\geq 4$ is an observation of Markl which was proved by Barnich
\cite{B}
( see the remark at the end of \cite{BFLS}).
A generalization of Markl's remark is available in a paper by Al-Ashhab
\cite{Samer}
and in that paper more explicit formulas are given for $l_1,l_2,l_3.$
Examination of these
formulas  provide the details needed for the calculations below.

First of all, we examine the mapping
$l_2: X_*\times X_*\longrightarrow X_*.$ Now
$ l_2 : A \times A $ $\longrightarrow X_*$
is  determined by $\tilde l_2: X_0\times X_0\longrightarrow X_0,$
consequently we
need only consider the restricted mapping:
\begin{eqnarray}
l_2: X_1\times X_0\longrightarrow X_1.
\end{eqnarray}
Moreover, since $X_0$ is a module over $k[[t]],$  $X_1$ is a module over
$k[[t]]t^2,$
and $\tilde l_2$ respects these structures
we need only consider its values on pairs $(a^{*}t^2,b)$ with
$a^{*}t^2\in X_1,b\in X_0.$   By  Theorem 2.2 of \cite{Samer}, we have
\begin{eqnarray}
l_2(a^{*}t^2,b)=-sl_2 l_1[(a^{*}t^2)\otimes b]\nonumber\\
=-sl_2[l_1(a^{*}t^2)\otimes b+(-1)^{\epsilon(a^*)}(a^{*}t^2)\otimes
l_1(b)]\nonumber\\
=-sl_2[(at^{2}\otimes b)]=-s[t^{2} l_{2}(a \otimes b)] \nonumber \\
=-s[t^2(\alpha_0(a,b)+\alpha_1(a,b)t)]\nonumber\\
=-s[\alpha_0(a,b)t^2+\alpha_1(a,b)t^3]\nonumber\\
=\alpha_0(a,b)^*t^2+\alpha_1(a,b)^*t^3\nonumber\\
=t^2(\alpha_0(a,b)^*+\alpha_1(a,b)^*t).
\end{eqnarray}

From this deduction, we  that the mapping $ l_{2} $ can essentially be
replaced  by
the modified map:
\begin{eqnarray}
\bar l_2: A[1] \times A\longrightarrow A[1][[t]],  \quad \quad \bar
l_2(a^*,b)=\alpha_{0}(a,b)^*+\alpha_{1}(a,b)^*t.
\end{eqnarray}
We clarify this remark below by showing that a new sh-Lie structure can
be obtained
with
$\bar l_2$ playing the role of $l_2.$

The next  mapping we examine is the mapping
\begin{eqnarray}
l_3: X_0 \times X_0\times X_0 \longrightarrow X_1
\end{eqnarray}
Since $l_{3} $ is $k[[t]]$-linear,we need only consider  mappings of the
type:
\newline $l_3: A\times A\times A \longrightarrow X_1$
where for $ x_1,x_2,x_3\in A,$
\begin{eqnarray}
l_3(x_1,x_2,x_3)=sl_2^2(x_1,x_2,x_3)\nonumber\\
=\sum_{\sigma \in unsh(2,1)}(-1)^{\sigma} sl_2(l_2(x_{\sigma(1)},
x_{\sigma(2)}),x_{\sigma(3)})\nonumber\\
=\sum_{\sigma \in unsh(2,1)}(-1)^{\sigma} sl_2
(\alpha_0(x_{\sigma(1)}, x_{\sigma(2)})+\alpha_1(x_{\sigma(1)},
x_{\sigma(2)})t,x_{\sigma(3)})\nonumber\\
  =\sum_{\sigma \in unsh(2,1)}(-1)^{\sigma}
s[\alpha_0(\alpha_0(x_{\sigma(1)},
x_{\sigma(2)}),x_{\sigma(3)})+\nonumber\\
   t\alpha_1(\alpha_0(x_{\sigma(1)}, x_{\sigma(2)}),x_{\sigma(3)})+
t\alpha_0(\alpha_1(x_{\sigma(1)}, x_{\sigma(2)}),x_{\sigma(3)})\nonumber\\
   +t^2 \alpha_1(\alpha_1(x_{\sigma(1)},
x_{\sigma(2)}),x_{\sigma(3)}),x_{\sigma(3)})\nonumber\\
  =s(\sum_{\sigma \in unsh(2,1)}(-1)^{\sigma}
\alpha_0(\alpha_0(x_{\sigma(1)},
x_{\sigma(2)}),x_{\sigma(3)})+\nonumber\\
+ t(\sum_{\sigma \in unsh(2,1)}(-1)^{\sigma}
\alpha_1(\alpha_0(x_{\sigma(1)},
x_{\sigma(2)}),x_{\sigma(3)})+\nonumber\\
\sum_{\sigma \in unsh(2,1)}(-1)^{\sigma} \alpha_0(\alpha_1(x_{\sigma(1)},
x_{\sigma(2)}),x_{\sigma(3)}))\nonumber\\
+t^2(\sum_{\sigma \in unsh(2,1)}(-1)^{\sigma}
\alpha_1(\alpha_1(x_{\sigma(1)},
x_{\sigma(2)}),x_{\sigma(3)}))]\nonumber\\
  =s((\alpha_0^2+t(\alpha_0\alpha_1+\alpha_1\alpha_0)+t^2
\alpha_1^2)(x_1,x_2,x_3))\nonumber\\
=s(t^2 \alpha_1^2 (x_1,x_2,x_3))\nonumber\\
  =-t^2 (\alpha_1^2 (x_1,x_2,x_3))^*
\end{eqnarray}
Or $l_3(x_1,x_2,x_3)=-t^2 ([\alpha_1,\alpha_1] (x_1,x_2,x_3))^*$ which
is precisely the
``first deformation obstruction class".}

Recall that we  know from Theorem 7 of \cite{BFLS} that we have an
sh-Lie structure.
The point of these calculations is that it enables us to obtain the
modified sh-Lie structure of Corollary 10 below and it is this structure
which is relevant to Lie algebra deformation.
Thus we already know that the mappings $l_1,l_2,l_3$  satisfy the relations:
\begin{eqnarray}
l_1l_2-l_1l_2=0   \\ {l_2}^2+l_1l_3+l_3l_1=0   \\{l_3}^2=0  \\
l_2l_3+l_3l_2=0.
\end{eqnarray}

Observe that if we let $\tilde X_*=\tilde X_1 \oplus \tilde X_0=A[1] [
[t] ] \oplus
A[[t]],$ then the formulas
defining $l_1,l_2, l_3$  defined on $X_*$ make sense on the new complex
$\tilde
X_*.$ Indeed
the calculations above show that  $l_1,l_3$ are uniquely determined by
their values
on ``constants" in the sense that they could be first defined on elements of
$A[1]\oplus A \subseteq A[1][[t]]\oplus A[[t]]$ and then extended to
$A[1][[t]]\oplus
A[[t]]$ using the fact that $l_1,l_3$ are required to be $k[[t]]$
linear. $l_2$ is not
obviously $k[[t]]$ linear. The whole point of corollary 10 below is that
the $sh$-Lie
structure defined by Theorem 10 can be redefined to obtain $sh$-Lie maps on
the graded space $\tilde X_*$ which are obviously
$k[[t]]$ linear and consequently this "new" structure is intimately
related to
deformation theory.
Thus, as we say above, the modified map
$\bar l_2$ can be extended to the new complex $\tilde X_*$ and is uniquely
determined by its values on ``constants".
If we denote the extensions of $l_1,l_3$ to $\tilde X_*$ by $\bar
l_1,\bar l_3,$
then clearly these
mappings satisfy the same relations (63)-(66) as the maps $l_1,l_2,l_3$ and
consequently
if we define $\bar l_n=0, n\geq 4$ it follows that  $(\tilde X_*,\tilde
l_1,\tilde
l_2,\tilde l_3,0,0 \cdots) $ is an $sh$-Lie algebra. This proves the
following
corollary.

\begin{corollary} There is an sh-Lie structure on $A[1][[t]]\oplus
A[[t]]$ whose
structure mappings
$\{\bar l_1,\bar l_2,\bar l_3, 0,\cdots\} $are precisely the mappings
$\{l_1,l_2,l_3, 0,\cdots\} $ when restricted to  $A[1][[t]]t^2\oplus
A[[t]].$ Moreover,
the structure mappings of $A[1][[t]]\oplus A[[t]]$ have the property
that they are
uniquely determined by their values on $A[1]\oplus A$ and $k[[t]]$
linearity.
\end{corollary}

From the  discussion above
the set of mappings  $\{\bar  l_1,\bar l_2,\bar l_3\}$ is essentially a
deformation
of an $sh$-Lie algebra. In addition, the construction    of
the mapping  $\bar l_2$ is  equivalent to defining an initial condition
for a  Lie
algebra deformation.

This means that a Lie algebra which can't be deformed in the category of
Lie algebra
may admit an $sh$-Lie algebra deformation by first imbedding it into an
appropriate
$sh$-Lie algebra.

\section{Conclusion}

Here, at  the end of the paper, we comment briefly  regarding the
developing role
chain extensions could play in mathematics and physics. We have seen
that chain
extensions may be regarded
as a generalization of ordinary algebraic extensions of Lie algebras.
This being the
case an interesting question is that of classifying such extensions
perhaps using
methods similar to  the  use of Chevalley-Eilenberg cohomology to
classify Lie
algebra extensions. More generally, it would be of interest
to develop a theory of chain extensions to parallel that of ordinary Lie
algebra
extensions and to consider applications of this theory beyond those
initiated in
this paper.

Another question of interest is that of the relationship between chain
extensions
and sh-Lie algebras. Preliminary calculations  suggest a link between
chain extensions and sh-Lie algebras and consequently that there may be
a way
of dealing with deformations of Lie algebras via chain extensions, but
more work is
required to be conclusive.
More generally, the notion of a chain extension has  provided us with a
new technique
for the investigation of deformation problems, but it's significance is
not yet
fully understood. Its role in physics is yet to be fully understood, but
we believe that
our reformulation of consistent deformations is an indication of its value.

\section*{Acknowledgments}
I would like to thank R.Fulp for useful  discussions throughout the writing
of this paper with special thanks
for his revision of the manuscript. I am also pleased to express my 
appreciation to
J.Stasheff who provided many thoughtful comments which led to substantial
improvement of  an early draft of my paper.

\newpage
%\bibliographystyle{amsplain}
%\bibliography{../Luminy/mybig}

\begin{thebibliography}{10}

\bibitem{Samer} S. AL-Ashhab,  {\em A class of strongly homotopy Lie
algebras with
simplied sh-Lie structures}, Los Alamos Archive,  math.RA/0308160

\bibitem{BBvD} F.A.Berends,G.J.H.Burgers and H.van Dam, Nucl.phys. B
{\bf 260}
(1985), 295

\bibitem{B} G.Barnich, {\em Brackets in the jet-bundle approach to field
theory},
Proceedings of the conference on Secondary Calculus and
Cohomological Physics, Contemp. Math. 219 (1998) 17-27.

\bibitem{BH} G.Barnich and M.Henneaux, Phys,Rev,Lett, {\bf 72} (1993), 1588

\bibitem{BFLS} G. Barnich,~R. Fulp,~T. Lada,~J. Stasheff,  {\em The sh
Lie structure
of Poisson brackets in field theory}, Commun. Math. Phys. {\bf
191}(1998),585-601

\bibitem{Brandt}F.Brandt, {\em Local BRST cohomology in the antifield
formalism},
unpublished lecture notes

\bibitem {G} M.Gerstenhaber,  {\em The cohomology structure of an
associate ring}
Annals of Mathematics {\bf 78} 267-288

\bibitem{HG} M.Hazewinkel and M.Gerstenharber, {\em Deformation theory
of Lie
algebras and structures and applications},
NATO Series C,Mathematics and Physics  Science {\bf 24}

\bibitem{HT} M.~Henneaux and C.~Teitelboim, {\em Quantization of {G}auge
{S}ystems},
Princeton Univ. Press, 1992.

\bibitem{KV} I.S. Krasilshchik and A.M. Vinogradov,  {\em Symmetries and
conservation laws for differential equations of mathematical physics},
Translations
of mathematical monographs {\bf 182,} AMS (1999)

\bibitem{ls}T.~Lada and J.D. Stasheff, {\em Introduction to sh {L}ie
algebras for
physicists}, Intern'l J. Theor. Phys. {\bf 32} (1993), 1087--1103.

\bibitem{Olver}P. Olver {\em Applications of Lie groups to differential
equations}
Graduate texts in Mathematics, Vol 107, Berlin-Heidelberg-New York:
Springer-Verlag,
1986

\bibitem{SS} M.Schlessinger,J.Stasheff,  {\em The {L}ie algebra
structure and
tangent cohomology and deformation}, J.Pure Appl. Algebra {\bf 89}
(1993), 231-235


\bibitem{s} J.Stasheff, {\em Deformation Theory and the
Batalin-Vilkovisky Master
Equation},
Proceedings of the Conference on Deformation Theory, etc. Ascona,
Switzerland (1996)
q-alg/9702012





\end{thebibliography}
%\input{jim218.bbl}

\ifx\undefined\bysame
\newcommand{\bysame}{\leavevmode\hbox to3em{\hrulefill}\,}
\fi

\end{document}